\documentclass[a4paper]{article}

\usepackage{amsfonts}
\usepackage{amsopn}
\usepackage{amsmath}
\usepackage{amsthm}
\usepackage{latexsym}
\usepackage{amssymb}
\usepackage{graphicx}
\usepackage{relsize}
\usepackage{subcaption}
\usepackage{lscape}

\usepackage[utf8]{inputenc}

\newcommand{\EE}{{\mathbb{E}}}

\newcommand{\RR}{{\mathbb{R}}}

\newcommand{\nP}{{\mathcal{P}}}

\newcommand{\re}{{\mathrm{e}}}

\newcommand{\sT}{{\mathsf{T}}}

\newcommand{\bxi}{{\boldsymbol{\xi}}}

\newcommand{\bx}{{\boldsymbol{x}}}

\newcommand{\by}{{\boldsymbol{y}}}

\newcommand{\sPot}{W}
\newcommand{\pPot}{V}
\newcommand{\sPotf}{w}
\newcommand{\pPotf}{v}

\renewcommand\d{\, \mathrm{d}}

\newcommand{\p}{\partial}

\renewcommand{\phi}{\varphi}

\newcommand{\ii}{{\rm i}}

\newcommand{\sklammer}[1]{\left\langle#1\right\rangle}
\newcommand{\norm}[1]{\left\|#1\right\|}
\newcommand{\dgl}[1]{\left\{\begin{aligned}
	#1
	\end{aligned}\right.}
\newcommand{\abs}[1]{\left|#1\right|}

\newcommand{\menge}[1]{\left\lbrace #1\right\rbrace }
\renewcommand{\vector}[1]{\begin{bmatrix} #1\end{bmatrix}}

\DeclareMathOperator{\Cov}{COV}

\DeclareMathOperator{\argmin}{argmin}

\newcommand{\vertiii}[1]{{\left\vert\kern-0.25ex\left\vert\kern-0.25ex\left\vert #1 \right\vert\kern-0.25ex\right\vert\kern-0.25ex\right\vert}}

\begin{document}

\title{Wave propagation in random media, parameter estimation and damage detection via stochastic Fourier integral operators}

\author{Michael Oberguggenberger\thanks{Unit of Engineering Mathematics, University of Innsbruck,
Technikerstra\ss e 13, 6020 Innsbruck,
Austria, (michael.oberguggenberger@uibk.ac.at)}
\and
Martin Schwarz\thanks{Unit of Engineering Mathematics, University of Innsbruck,
Technikerstra\ss e 13, 6020 Innsbruck, Austria,
(math@mschwarz.eu)}
}

\date{}
\maketitle

\begin{abstract}
This paper presents a new approach to modelling wave propagation in random, linearly elastic materials, namely by means of Fourier integral operators (FIOs). The FIO representation of the solution to the equations of motion can be used to identify the elastic parameters of the underlying media, as well as their statistical hyperparameters in the randomly perturbed case. A stochastic version of the FIO representation can be used for damage detection. Hypothesis tests are proposed and validated, which are capable of distinguishing between an undamaged and a damaged material, even in the presence of random material parameters.

The paper presents both the theoretical fundamentals as well as a numerical experiment, in which the applicability of the proposed method is demonstrated.
\end{abstract}


\section{Introduction}

Modelling wave propagation in elastic solids goes back to the first half of the 19th century. According to the short historical survey in \cite{Achenbach1976}, the  work of Raleigh, Lamb and Love in the first half of the 20th century sparked a huge development of the subject first in seismology \cite{Bleistein2001} and then in the material sciences \cite{ewins1984modal} and geotechnics \cite{Ostrovsky:2001}. In the first case, elastic waves, often sound waves, are used for detecting the internal structure of the earth crust; in the second case, ultrasonic measurements are used for damage detection, as well as mechanically excited waves. Under the heading of vibration based condition monitoring  \cite{Carden2004}, a wide range of methods has been developed, from modal and spectral analysis to finite element based dynamical models and more. The reader is referred e.g. to the monographs \cite{brandt2011noise,Fahy2007,gopalakrishnan2007spectral,Stepinski2013}
and the survey papers \cite{Carden2004,Deraemaeker2010,Doebling1998,Fan2011,Farrar:2001,Marcantonio2019,Ostachowicz2010,Palacz2018}.

In the past two decades, additional attention has been given to the random structure of soils and materials and the modelling of wave propagation in randomly layered or stochastically perturbed media; a sample of references from seismology, geotechnics, and material sciences is \cite{Borcea:2016,Demmie2016,Fouque2007,Manolis2002,NairWhite:1991}.

This paper proposes a new approach to modelling wave propagation in random elastic materials, based on Fourier integral operators (FIOs). This type of operators has been introduced by \cite{Hoermander:1971} in the 1970s and is now widely used for representing solutions to hyperbolic partial differential equations. The operators are of the form
\begin{equation}\label{eq:FIO}
u(\bx) = A[f](\bx)=\int\int \re^{\ii\Phi(\bx,\by,\bxi)} a(\bx,\by,\bxi) f(\bx) \d \by\d\bxi,
\end{equation}
where $\Phi$ is the \emph{phase function}, which encodes the geometry of wave propagation, and $a$ is the \emph{amplitude function}, representing the frequency content. Typically, $f(\bx)$ is the source and $u(\bx)$ the response of the medium; time dependence is not shown explicitly in formula \eqref{eq:FIO}.
A further novelty of the presented approach is that randomness of the material will be modelled through the stochastic phase functions and amplitudes of the representing Fourier integral operators. This is in contrast to the common approach in which randomness in the material is incorporated by taking the coefficients of the equations of motion as random fields \cite{Matthies:2008,Mishra:2016}. The proposed alternative approach has the advantage that the stochastic system response, described through \eqref{eq:FIO}, can be analyzed and simulated in a rather explicit form. Fourier integral operators will be used for
\smallskip

(a) propagating acoustic waves through the material and computing the dynamic response to localized excitations;
\smallskip

(b) identifying the elastic parameters of the material (Young's modulus, Poisson's ratio) as well as their statistical hyperparameters (variances and correlation lengths);
\smallskip

(c) designing statistical tests detecting significant changes in the parameters (due to damage or fatigue).
\smallskip

In short, the approach proceeds as follows. The general assumption is that the material is linearly elastic with inherent randomness. The spatial randomness of the material parameters is caused by small random perturbations of otherwise constant nominal values. The first step is to set up the Fourier integral operators solving the equations of motion with constant material parameters. These parameters appear explicitly in the phase functions and amplitudes of the operators and -- in the second step -- are replaced by spatial random fields there. The material parameters and their statistical hyperparameters are then estimated by comparing the FIO solution to measured data in a finite number of locations (sensors). This results in a stochastic Fourier integral operator model of the propagation of waves in the given elastic material, interpreted as the undamaged state. In a third step, this model can be used to generate a large number of system responses of the undamaged state by Monte Carlo simulation. This sample can be interpreted as a realization of the null hypothesis of undamaged material, against which later measurements of the possibly damaged material are tested.

Statistical methods for damage detection have been proposed in the literature before \cite{Farrar:2010,Nichols:2016}, also in combination with Monte Carlo sampling \cite{Zhao2019}.
Compared to other numerical methods (finite differences, finite elements, modal analysis) the advantage of the FIO based approach lies in the following. Usually, in parameter identification and damage detection one has only a few sensors at certain locations of the material. These sensors measure the dynamic system response, which is a time-dependent signal. Thus one does not have complete spatial information at hand. This setup motivates Fourier integral operators, since Fourier integral operators can simulate the response to an excitation of an elastic material at single points without computing the solution on the whole domain. This makes parameter estimation very efficient. Only setting up and calibrating the baseline model of undamaged material requires an optimization procedure and a larger amount of model evaluations. Generating Monte Carlo samples of the baseline response and comparing it with (features of) measured data is computationally inexpensive.

This paper serves to set up the theoretical concepts and to demonstrate the applicability in the case of isotropic, linearly elastic materials. In this case, the Helmholtz decomposition allows one to decouple the equations of motion into four scalar wave equations for the wave potentials, which in turn are conveniently solved by Fourier integral operators. However, the approach works for other materials (e.g. orthotropic materials) as well. The hypothesis tests will be validated by artificially generated measured data, computed by a finite element model.

The computational advantages of the FIO method come with certain limitations on the range of applicability. In fact, FIOs of the form \eqref{eq:FIO} can model forward propagation of stochastically perturbed waves, but the representation does not take into account internal reflections due to local random changes of the material properties. Thus local scattering or mode conversions are not modelled. Consequently, the presented perturbative approach is valid for random media in the weakly heterogeneous scaling regime only. (It is quite common in the theory of wave propagation in random media that different regimes are analyzed by means of different models, see e.g. \cite{Fouque2007,NairWhite:1991,Wu:1988}.)

The structure of the paper is as follows. The second section serves to present the theoretical background, the Helmholtz decomposition and the solution of scalar wave equations by means of Fourier integral operators in three space dimensions. The concept of stochastic phase and amplitude functions will be briefly introduced. Detailed formulas and explicit expressions for the representing Fourier integral operators will be given in the case of a plane strain problem.

In the third section, it is assumed that time-dependent measured data at certain locations are available, and it is shown how to estimate the material parameters and their statistical hyperparameters and how to calibrate the baseline stochastic FIO-model.

The fourth section briefly describes how statistical hypothesis tests can be designed with the help of the baseline stochastic FIO-model.

The last section sets out to validate the proposed procedure in a numerical example. It is assumed that a block of aluminum is excited by a harmonic line force. The laboratory measurements are represented by a finite element simulation with stochastically perturbed Lam\'e parameters. It is shown that the FIO-solution is able to reproduce the results of the FE-simulation for unperturbed constant Lam\'e parameters accurately. Furthermore, if the material parameters in the FE-simulation are randomly perturbed, it is demonstrated how the nominal values, standard deviations and correlation lengths of the Lam\'e parameters can be estimated. This results in the calibrated stochastic FIO-model of the undamaged material. The stochastic FIO-model of the undamaged state is used to generate a Monte Carlo sample of possible system responses of the undamaged material. Four scenarios of material states are considered: undamaged material, stiff material, soft material, or presence of a crack. The latter three are considered as degradation, fatigue, or damage. Three hypothesis tests are applied, and it is shown that they can distinguish between all four scenarios. Finally, the performance of the tests is evaluated with respect to false classification rates. Repeating the tests with different samples of the underlying random fields, it is shown that the tests behave as designed.

A rough quantification of the range of applicability, following the classification of scaling regimes from \cite{Fouque2007}, concludes this section.

The paper ends with the conclusion, followed by an appendix containing additional graphical and tabular illustrations of the results referred to in the text.

The paper is based on the work \cite{Schwarz:Thesis:2019} of the second author. The mathematical background of stochastic Fourier integral operators has been published in \cite{MOMS:2020}. The paper is a continuation of \cite{MOMS:2019}, where the propagation of ultrasonic waves in one-dimensional media was addressed.
The authors would like to point out that another aspect of Fourier integral operators combined with stochastic processes, even for nonlinear equations, has been worked out in a series of papers \cite{Ascanelli:2018,Coriasco:2019,Coriasco:2020}. However, in these papers the phase functions are taken deterministic, the stochastic processes appear in the driving forces.

\section{Fourier integral operators in linear elasticity}
\label{sec:FIO_solution_3d_linelst}

\subsection{Helmholtz decomposition}
\label{subsec:Helmholtz}

This subsection serves to motivate the occurrence of FIOs in the solution of the equations of motion of linear elasticity. The displacement $\boldsymbol{u}$ of a linear elastic, homogeneous, isotropic body is governed by the equations of motion
\begin{align}
\label{eqn:3dlinelst}
\dgl{\partial_{tt}\boldsymbol{u}(\boldsymbol{x},t) -\frac{\lambda+\mu}{\rho}(\nabla\nabla\cdot\boldsymbol{u}(\boldsymbol{x},t))-\frac{\mu}{\rho} \Delta \boldsymbol{u}(\boldsymbol{x},t) =\boldsymbol{f}(\boldsymbol{x},t)\\
	\boldsymbol{u}(\bx,0) = \boldsymbol{0} , \quad \p_t \boldsymbol{u}(\bx,0) = \boldsymbol{0} }
\end{align}
where $\boldsymbol{f}$ is the body force, $\lambda$ and $\mu$ the Lam\'e constants and $\rho$ the density; $\Delta$ denotes the Laplace operator. The classical Helmholtz decomposition allows one to decompose this system of equations into four decoupled wave equations of the form
\[
\dgl{\p_{tt} U(\bx,t)-c^2\Delta U(\bx,t)= F(\bx,t)\\U(\bx,0)=0, \quad \p_t U(\bx,0)=0}
\]
where $c$ is either the lateral wave speed ($c^2=(\lambda+2\mu)/\rho$) or transversal wave speed ($c^2=\mu/\rho$), and $u$ denotes the pressure potential or one of the components of the shear potential. Details about the decomposition can be found, e.g., in \cite{Achenbach1976}; the case of two space dimensions will be elaborated in Subsection~\ref{subsec:solving_2D_model}. For simplicity, assume that $F$ can be written as $F(\bx,t)= g(\bx) h(t)$. Taking the Fourier transform with respect to $\bx$,
\[
  \widehat{U}(\bxi,t) = \int_{\RR^3}\re^{-\ii\sklammer{\bx,\bxi}}U(\bx,t) \d\bx
\]
yields the ordinary differential equation
\begin{align*}
\dgl{\p_{tt} \widehat{U}(\bxi,t) + c^2 \norm{\bxi}^2 \widehat{U}(\bxi,t) = \widehat g(\bxi) h(t) \\
	\widehat U(\bxi,0)=0, \quad \p_t \widehat U(\bxi,0)=0}
\end{align*}
which is solved by
\begin{align}
\widehat{U}(\bxi,t)= \frac{ \widehat{g}(\bxi)}{2 \ii c \norm{\bxi}} \bigg( \int_0^t h(s)\re^{\ii c \norm{\bxi} (t-s)}  \d s - \int_0^t h(s) \re^{-\ii c \norm{\bxi} (t-s)}  \d s \bigg).
\end{align}
Here $\sklammer{\bx,\bxi}$ denotes the Euclidean scalar product and $\norm{\bxi}$ the Euclidean norm.
Applying the inverse Fourier transform yields the solution operator, represented by means of oscillatory integrals
\begin{align}\label{eqn:solution_operator_FIO}
\begin{aligned}
U(\bx,t)= \frac{1}{(2\pi)^3}\bigg(& \int_{\RR^3} \int_0^t \re^{\ii\sklammer{\bx,\bxi}+ \ii  c \norm{\bxi} (t-s)} \frac{  \widehat{g}(\bxi) h(s)}{2 \ii c \norm{\bxi}} \d s \d \bxi\\
-&\int_{\RR^3} \int_0^t \re^{\ii\sklammer{\bx,\bxi}- \ii  c \norm{\bxi} (t-s)} \frac{  \widehat{g}(\bxi) h(s)}{2 \ii c \norm{\bxi}} \d s \d \bxi\bigg),
\end{aligned}
\end{align}
These integrals can be seen as time-integrated Fourier integral operators. For example, the first integral contains
\begin{equation}\label{eq:FIOterm}
  \int_{\RR^3}\re^{\ii\sklammer{\bx,\bxi}+ \ii  c \norm{\bxi} (t-s)} \frac{\widehat{g}(\bxi)}{2 \ii c \norm{\bxi}} \d \bxi
  = \int_{\RR^3}\int_{\RR^3}\re^{\ii\sklammer{\bx-\by,\bxi}+ \ii  c \norm{\bxi} (t-s)} \frac{g(\by)}{2 \ii c \norm{\bxi}}\d\by \d \bxi
\end{equation}
which is of the form \eqref{eq:FIO} with $\Phi(\bx,\by,\bxi) = \sklammer{\bx-\by,\bxi}+  c \norm{\bxi} (t-s)$
and with $a(\bx,\by,\bxi) = (2 \ii c \norm{\bxi})^{-1}$.
In general, this integral does not converge as it stands, but has to be interpreted as an \emph{oscillatory integral}, as e.g. described in \cite{Shubin:01}.
(Actually, an additional regularization at $\bxi = {\mathbf 0}$ is required to fully conform with the theory of Fourier integral operators.)
However, if $g$ is sufficiently regular so that $\widehat{g}$ is absolutely integrable then \eqref{eq:FIOterm} can be directly evaluated by iterated integration.

\subsection{Stochastic perturbations} \label{subsec:stoch}
In the derivation of the solution operator in the previous subsection, the material parameters were assumed to be constant. Generally, even in an undamaged material,
the material parameters show spatially random deviations. (See e.g. the ultrasonic measurements of composite plates in \cite{MOMS:2019}.) As outlined in the introduction, the random structure of the material will be accounted for by adding random perturbations of the phase function and the amplitude function in the solution operator. These random perturbations take the form of spatial random fields, to be calibrated by measurements.

This has the advantage that evaluating the solution by means of stochastically perturbed Fourier integral operators is much simpler than modelling the coefficients in system
\eqref{eqn:3dlinelst} and then computing the stochastic solution by standard numerical methods.

To construct the perturbed Fourier integral operator, the following strategy is used. Let
\[
\iint \re^{\ii\Phi(\boldsymbol{x},\boldsymbol{\xi},t)-\ii\boldsymbol{\xi}\cdot\boldsymbol{y}} a(\boldsymbol{x},\boldsymbol{\xi}) g(\boldsymbol{y}) \mathrm{d} \boldsymbol{y}\d\boldsymbol{\xi}
\]
be a term of the solution operator for the wave equation with constant wave speed, as in \eqref{eq:FIOterm}. Adding random perturbations $r_1(\boldsymbol{x},\boldsymbol{\xi},t)$, $r_2(\boldsymbol{x},\boldsymbol{\xi})$ results in the stochastic Fourier integral operator
\begin{equation}
\iint \re^{\ii(\Phi(\boldsymbol{x},\boldsymbol{\xi},t)+r_1(\boldsymbol{x},\boldsymbol{\xi},t))-\ii\boldsymbol{\xi}\cdot\boldsymbol{y}} (a(\boldsymbol{x},\boldsymbol{\xi})+r_2(\boldsymbol{x},\boldsymbol{\xi})) g(\boldsymbol{y}) \mathrm{d} \boldsymbol{y}\d\boldsymbol{\xi} \label{eqn:stoch_pert_fio_waves}.
\end{equation}
This representation has a physical interpretation, namely $r_1$ describes the perturbation geometry of the propagation, while $r_2$ is the perturbation amplitude of the wave number.

If the input $g$ is not regular enough, for example, a rectangular pulse or even a Dirac delta function, one needs to ensure that the perturbed integral still defines an oscillatory integral. Thus, one has to make sure that $r_1$ and $r_2$ are of the right class, i.e., $\Phi+r_1$ must lie in an appropriate space of phase functions (see \cite{MOMS:2020}) and $a+r_2$ must belong to the H\"ormander symbol class $S_{\varrho,\delta}^m$ with $0<\varrho\leq1, 0\leq \delta <1$ (see e.g. \cite{Shubin:01}).
This difficulty will be avoided in the sequel by assuming that $g(\by)$ in \eqref{eqn:stoch_pert_fio_waves} is sufficiently smooth and has a sufficiently rapid decay rate at infinity, so that the Fourier integral operators are given by converging iterated integrals, which can be directly numerically evaluated.

\subsection{A plane strain model problem} \label{subsec:solving_2D_model}

For the sake of simplicity in presentation and computation, the methods of the previous subsections will be elaborated in detail in the case of two space dimensions.
Consider a linearly elastic, homogeneous, isotropic body occupying the whole three-dimensional space.
Assume that the force term $\boldsymbol f$ in \eqref{eqn:3dlinelst} is constant along the $z$-axis, i.e., depends only on $x$, $y$ and $t$.
Since the body is assumed to be at rest at $t=0$, the solution $\boldsymbol u$ to \eqref{eqn:3dlinelst} is constant along the $z$-axis as well, and it is of the form
\[
\boldsymbol{u}(x,y,z,t)=[u_1(x,y,t), u_2(x,y,t),0]^\sT.
\]
Without loss of generality, one may assume that the force is applied only in $y$-direction. More formally, it is assumed that the force is of form
\[
\boldsymbol{f}(x,y,z,t)=[0,g(x,y)h(t),0]^\sT,
\]
where $g$ is the support of the force, and $h$ is the time-dependent force term.
The problem is reduced to two dimensions, where the displacements $u_1(x,y,t)$, $u_2(x,y,t)$ satisfy the equations
\begin{align} \label{eqn:linelst2d_full_eqn}
\p_{tt} \vector{u_1\\u_2}-\frac{1}{\rho}\begin{bmatrix}
({\lambda+2\mu}) \p_x^2 + {\mu}{}\p_y^2 & (\lambda+\mu) \p_{xy}\\
(\lambda+\mu) \p_{xy} & ({\lambda+2\mu}) \p_y^2 + {\mu}{}\p_x^2
\end{bmatrix} \vector{u_1\\u_2}  = \vector{0\\g h}.
\end{align}
In this case, the Helmholtz decomposition is performed as follows, see e.g. \cite{Achenbach1976}. There is a pressure  potential $V$ and a shear potential $W$ such that
\begin{equation}\label{eq:u1u2}
\begin{bmatrix}
u_1(y,x,t)\\u_2(x,y,t)
\end{bmatrix} = \nabla \pPot (x,y,t) + \begin{bmatrix}
\partial_y \sPot(x,y,t)\\-\partial_x \sPot(x,y,t)
\end{bmatrix}.
\end{equation}
The potentials satisfy the wave equations
\begin{align}
\begin{aligned}\label{eqn:wave2dpot}
\partial_{tt} \pPot(x,y,t) - c_l^2 \Delta \pPot(x,y,t) =\pPotf(x,y)h(t), \\ \partial_{tt} \sPot(x,y,t) - c_s^2 \Delta\sPot(x,y,t) =\sPotf(x,y)h(t).
\end{aligned}
\end{align}
with $c_l^2= (\lambda+2\mu)/\rho$ and $c_s^2=\mu/\rho$, where the functions $\pPotf$ and $\sPotf$ are to be determined from
\[ \begin{bmatrix}
0\\g(x,y)
\end{bmatrix} = \nabla \pPotf (x,y) + \begin{bmatrix}
\partial_y \sPotf(x,y)\\-\partial_x \sPotf(x,y)
\end{bmatrix}.\]
Making the ansatz $v(x,y)=\p_y \varphi(x,y)$ and $w(x,y)=-\p_x \varphi(x,y)$ for some function $\phi(x,y)$, one has that
\[ \nabla \pPotf (x,y) + \begin{bmatrix}
\partial_y \sPotf(x,y)\\-\partial_x \sPotf(x,y)
\end{bmatrix}= \begin{bmatrix}
0\\\Delta \phi(x,y)
\end{bmatrix}.\]
Thus, by solving the two-dimensional Poisson equation
\begin{align}\label{eqn:inverting_laplacean}
 \Delta \phi = g,
\end{align}
with the growth condition $\lim_{\norm{(x,y)}\to\infty} \phi(x,y)=0$, one uniquely obtains $\phi$ and thus $v$ and $w$.

In more detail, if $g(x,y)$ is the Dirac delta function $\delta(x,y)$, then $\phi$ is the fundamental solution of the two-dimensional Laplace operator,
\[
\phi(x,y)=(2\pi)^{-1} \log(\sqrt{x^2+y^2}).
\]
For general compactly supported $g(x,y)$, the solution is obtained by convolution with the fundamental solution.

To continue solving \eqref{eqn:linelst2d_full_eqn}, the spatial Fourier transform is applied in equation \eqref{eqn:wave2dpot}, similarly to Subsection~\ref{subsec:Helmholtz}. This leads to the ordinary differential equations
\begin{align*}
\begin{aligned}
\partial_{tt} \widehat{\pPot}(\xi,\eta,t) = -c_l^2 (\xi^2+\eta^2) \widehat{\pPot}(\xi,\eta,t) +\widehat{\pPotf}(\xi,\eta)h(t) \\
\partial_{tt} \widehat{\sPot}(\xi,\eta,t) = -c_s^2 (\xi^2+\eta^2) \widehat{\sPot}(\xi,\eta,t) +\widehat{\sPotf}(\xi,\eta)h(t).
\end{aligned}
\end{align*}
Solving the ordinary differential equations and applying the inverse Fourier transform one ends up with the following solution operators:
\begin{align}
\begin{aligned}\label{eqn:2dphi}
\pPot(x,y,t)=&\frac{1}{(2\pi)^2}\int_0^t \int_{\RR^2} \re^{\ii x \xi+ \ii y \eta}\, \frac{\sin\big( c_l(t-s)\sqrt{\xi^2+\eta^2}\big)}{c_l\sqrt{\xi^2+\eta^2}}\,{ \widehat{\pPotf}(\xi,\eta)}h(s) \d (\xi,\eta) \d s
\end{aligned}\\
\intertext{and}
\begin{aligned} \label{eqn:2dpsi}
\sPot(x,y,t)=&\frac{1}{(2\pi)^2}\int_0^t \int_{\RR^2} \re^{\ii x \xi+ \ii y \eta}\, \frac{\sin\big( c_s(t-s)\sqrt{\xi^2+\eta^2}\big)}{c_s\sqrt{\xi^2+\eta^2}}\,{ \widehat{\sPotf}(\xi,\eta)}h(s) \d (\xi,\eta) \d s.
\end{aligned}
\end{align}
Instead of inverting the Laplace operator as in \eqref{eqn:inverting_laplacean}, one can simply assert that
\[
 \widehat{\pPotf}(\xi,\eta)=\frac{\ii \eta}{-(\xi^2+\eta^2)} \widehat{g}(\xi,\eta),\quad \widehat{\sPotf}(\xi,\eta)=\frac{-\ii \xi}{-(\xi^2+\eta^2)} \widehat{g}(\xi,\eta).
\]
Together with the relation
\[ \begin{bmatrix}
\widehat{u}_1(\xi,\eta)\\
\widehat{u}_2(\xi,\eta)
\end{bmatrix} = \begin{bmatrix}
\ii \xi& \ii \eta\\ \ii\eta &-\ii\xi
\end{bmatrix} \ \begin{bmatrix}
\widehat{\pPot}(\xi,\eta)\\\widehat{\sPot}(\xi,\eta)
\end{bmatrix}\]
from \eqref{eq:u1u2}, the final solution is obtained as
\begin{align} \label{eqn:linelst2d_num_sol_u}
u_{{1,2}}(x,y,t)= \frac{1}{(2\pi)^2}\int_{\RR^2} \re^{\ii x \xi+ \ii y \eta}  \widehat{u}_{{1,2}} (\xi,\eta,t) \d(\xi,\eta),
\end{align}
where
\begin{align}
\begin{aligned} \label{eqn:linelst2d_num_sol_u1_hat}
\widehat{u}_1(\xi,\eta,t)&=\int_0^t\left(\frac{\sin\big( c_l(t-s)\sqrt{\xi^2+\eta^2}\big)}{c_l\sqrt{\xi^2+\eta^2}} \frac{ \xi\eta \ \widehat g(\xi,\eta)}{\xi^2+\eta^2}h(s)  \right. \\
&\qquad \ - \left. \frac{\sin\big( c_s(t-s)\sqrt{\xi^2+\eta^2} \big)}{c_s\sqrt{\xi^2+\eta^2}} \frac{ \xi \eta \ \widehat g(\xi,\eta)}{\xi^2+\eta^2}h(s)  \right) \d s
\end{aligned}
\end{align}
and
\begin{align}
\begin{aligned} \label{eqn:linelst2d_num_sol_u2_hat}
\widehat{u}_2(\xi,\eta,t)&=\int_0^t\left(  \frac{\sin\big( c_l(t-s)\sqrt{\xi^2+\eta^2}\big)}{c_l\sqrt{\xi^2+\eta^2}} \frac{ \eta^2 \ \widehat g(\xi,\eta)}{\xi^2+\eta^2}h(s) \right. \\
&\qquad \ + \left.  \frac{\sin\big( c_s(t-s)\sqrt{\xi^2+\eta^2}\big)}{c_s\sqrt{\xi^2+\eta^2}} \frac{ \xi^2 \ \widehat g(\xi,\eta)}{\xi^2+\eta^2}h(s) \right) \d s.
\end{aligned}
\end{align}
These integrals are convergent, since there are no poles left anymore.

\emph{Remark.}	The solution operator \eqref{eqn:linelst2d_num_sol_u} can be still interpreted as a sum of Fourier integral operators. Indeed, letting
$\alpha=c_{l,s} (t-s) \sqrt{\xi^2+\eta^2}$ and observing that $\sin(\alpha)=(\re^{\ii\alpha}-\re^{-\ii\alpha})/2\ii$, the respective phase functions are seen to be of the form
\[
\begin{aligned}
\Phi_l(x,y,\xi,\eta,s,t) = x\xi+y\eta\pm c_l(t-s)\sqrt{\xi^2+\eta^2},\\
\Phi_s(x,y,\xi,\eta,s,t) = x\xi+y\eta\pm c_s(t-s)\sqrt{\xi^2+\eta^2}.\\
\end{aligned}
\]
Before describing the shape of the stochastic perturbations of the phase functions, and for later reference,  a switch of notation from the Lam\'e parameters $\lambda$ and $\mu$
to the more commonly used elastic material parameters, Young's modulus $E$ and Poisson's ratio $\nu$, is convenient. The parameters are related by
\[
\lambda=\frac{E\nu}{(1+\nu)(1-2\nu)}, \qquad \mu=\frac{E}{2 (\nu+1)},
\]
respectively. The pressure and shear wave speeds are expressed as
\begin{equation}\label{eq:wavespeeds}
c_l=\sqrt {\frac{ E}{\rho}{\frac { \left( 1-\nu \right)}{  \left( 1+\nu \right)  \left(1-2\,\nu \right)}}},\qquad c_s=\sqrt{\frac{E}{\rho} \frac{1}{2(\nu+1)}},
\end{equation}
respectively. The stochastic perturbations of the phase functions are obtained by making $c_l$ and $c_s$ space dependent and setting
\begin{equation}\label{eq:randwavespeeds}
\begin{aligned}
c_l(x,y)&=\sqrt {\frac{1}{\rho}{\frac { \left( 1-\nu(x,y) \right) E(x,y)}{ \left( 1+\nu(x,y) \right)  \left(1-2\,\nu(x,y) \right)}}}, \\
c_s(x,y)&=\sqrt{\frac{1}{\rho} \frac{E(x,y)}{2(\nu(x,y)+1)}},
\end{aligned}
\end{equation}
that is, the constant parameters $E$ and $\nu$ from formula \eqref{eq:wavespeeds} are replaced by random fields $E(x,y)$ and $\nu(x,y)$ of the form
\begin{equation}\label{eq:RFEnu}
E(x,y)=E_{\text{mean}}+R_E(x,y), \qquad
\nu(x,y)=\nu_\text{mean}+R_\nu(x,y).
\end{equation}
Here $R_E$ and $R_\nu$ are homogeneous, mean zero random fields, whose statistical parameters will later have to be calibrated by measurements. Note that the density $\rho$ enters in $c_l$ and $c_s$ only through the quotient $E/\rho$. Thus, for the purpose of this analysis, it is no restriction of generality to take $\rho$ constant, putting all randomness formally into $E$ and $\nu$.

Strictly speaking, the insertion of non-constant parameters in \eqref{eq:randwavespeeds} violates the derivation of the underlying equations in \eqref{eqn:linelst2d_full_eqn}; otherwise occurring spatial derivatives of the parameters are neglected. However, this approximation is in line with the general limitation of the method to small random perturbations, as detailed in Subsection~\ref{subsec:advantages}.

\emph{Remark.} In reality, elastic bodies occupy a finite region. Then system \eqref{eqn:3dlinelst} will have to be complemented by boundary conditions. However, in the intended applications, the force is localized in a small spatial region. Thus, for short time, wave propagation can be modelled by the full-space system, as long as the waves do not reach the boundary.
(For a discussion of distances and travelling times in the numerical example of Section~\ref{sec:numerical} see Subsection~\ref{subsec:advantages}.)

\section{Stochastic parameter estimation}
\label{sec:Parameter_estimation}

The random fields \eqref{eq:RFEnu} associated with the modulus of elasticity and Poisson's ratio are determined by their type of distribution and their statistical parameters. This section presents the proposed procedure how to identify these parameters using the Fourier integral operator representation of the solution.

It will be assumed that the random fields $R_E$ and $R_\nu$ are independent, homogeneous, centered and Gaussian with autocovariance functions
\begin{align}
\begin{aligned}
\label{eq:basicautocov}
\Cov\big(R_E(x_1,y_1),R_E(x_2,y_2)\big)&=\sigma^2_E\, \re^{-r_{12}/L_E}, \\ \Cov\big(R_\nu(x_1,y_1),R_\nu(x_2,y_2)\big)&=\sigma^2_\nu\, \re^{-r_{12}/L_\nu}
\end{aligned}
\end{align}
where $r_{12} = \sqrt{(x_1-x_2)^2+(y_1-y_2)^2}$ is the Euclidean distance. Here $\sigma_E $, $\sigma_\nu$ and $L_E$, $L_\nu$ denote the respective standard deviations and correlation lengths. However, these assumptions are not obligatory. Different types of autocovariance functions can easily be implemented, and the same methods are also applicable for non-Gaussian random fields.

The randomness of the material is thus described by the six parameters $E_{\text{mean}},\sigma_E, L_E, \nu_{\text{mean}}, \sigma_\nu, L_\nu$.
The model scenario is that the time-dependent system output is recorded at $n$ sensor locations. In practice, this will be a record of measured data. For the present study, the record is artificially generated by computer simulation. In any case, it is assumed that the displacements $u_1$ and $u_2$ are known at the locations $(x_j,y_j),j=1,\ldots,  n$ for all times $t$ in the time interval under consideration. The task is to estimate the six parameters by comparing the data with the solution obtained by the Fourier integral operator representation \eqref{eqn:linelst2d_num_sol_u}.

\subsection{Estimating the nominal values}
\label{subsec:nomval}

The first task is to estimate the nominal values $E_\text{mean} $ and $\nu_\text{mean}$ based on the data given in the sensor locations. This will be done by fitting FIO solutions \eqref{eqn:linelst2d_num_sol_u} with constant parameters $E_0,\nu_0$ in \eqref{eq:wavespeeds}.

To shorten notation, let $\boldsymbol{g}_{j}(t)=\boldsymbol{u}_{\text{meas}}(x_j,y_j,t)$ be the measured signal in $(x_j,y_j)$, recorded over the time interval
$[0, T_{\text{max}}]$. Thus  $\boldsymbol{g}=[g_{d,j}], d=1,2$ is a matrix-valued function. The FIO solution $\boldsymbol{h}_{j}(t,E_0,\nu_0)=\boldsymbol{u}_{\text{FIO}}(x_j,y_j,t,E_0,\nu_0)$ solves the deterministic, constant coefficient case with parameters $E_0$ and $\nu_0$, and $\boldsymbol{h}=[h_{d,j}],d=1,2$.

Goodness of fit is measured in terms of the weighted norm
\begin{align}\label{eqn:error_function_parameest}
\vertiii{\boldsymbol{f}} = \sum_{j=1}^n \norm{\boldsymbol{f}_{j}}_{2}=\sum_{j=1}^{n} \sqrt{\sum_{d=1}^2w_d\int_0^{T_{\text{max}}} f_{d,j}^2(t) \d t},
\end{align}
and $E_0$, $\nu_0$ are estimated as solutions to the minimization problem
\begin{equation}\label{eq:argmin}
[E_{0,\text{opt}},\nu_{0,\text{opt}}]=\argmin_{{E}_0,{\nu}_0}\vertiii{\boldsymbol{g}-\boldsymbol{h}(\cdot,{E}_0,{\nu}_0)}.
\end{equation}
The weights $w_1$ and $w_2$ are chosen as the maximum of the measured displacements in $x$-direction and in $y$-direction, respectively, at each sensor location and over the whole time interval $[0, T_{\text{max}}]$. This ensures that both directions enter with equal importance in the error norm.

If a sample of size $N$ of the measured signals $\boldsymbol{u}_{\text{meas}}(x_j,y_j,t)$ is available, one can increase the robustness of the algorithm by repeating the calibration $N$ times, resulting in the estimates $E_{0,\text{opt}}^{(1)}, \ldots, E_{0,\text{opt}}^{(N)}$ and
$\nu_{0,\text{opt}}^{(1)}, \ldots, \nu_{0,\text{opt}}^{(N)}$, and setting:
\begin{equation}\label{eq:EstarNustar}
 E_\text{mean}^*= \frac{1}{N} \sum_{k=1}^N E_{0,\text{opt}}^{(k)},\qquad
  \nu_\text{mean}^*= \frac{1}{N} \sum_{k=1}^N \nu_{0,\text{opt}}^{(k)}.
\end{equation}
Numerical experiments showed that $E_\text{mean}^*$ and $\nu_\text{mean}^*$ are good estimators of $E_\text{mean}$ and $\nu_\text{mean}$, respectively, even with small sample size $N$ (see Subsection\;\ref{subsec:parest}).

\subsection{Estimating the standard deviations}
\label{subsec:stds}

The next task is to estimate the standard deviations $\sigma_E$ and $\sigma_\nu$. To obtain the required data, the nominal values are first estimated for each sensor location separately:
\begin{align*}
[E_{0,\text{opt},j},\nu_{0,\text{opt},j}]=\argmin_{{E}_0,{\nu}_0}\norm{\boldsymbol{g}_j(\cdot)-\boldsymbol{h}_j(\cdot,{E}_0,{\nu}_0)}_{2}
\end{align*}
for $j = 1,\ldots, n$.
The estimator of the  variance of $E_{0,\text{opt}}$ and $\nu_{0,\text{opt}}$ will be denoted by $\sigma_{E_{0,\text{opt}}}^2$ and is computed as follows:
\begin{equation}\label{eq:sigEopt}
\sigma_{E_{0,\text{opt}}}^2 = \frac{1}{N\, n-1}  \sum_{j=1}^n \sum_{k=1}^N \left (E_{0,\text{opt},j}^{(k)} - E_\text{mean}^* \right)^2,
\end{equation}
invoking the hypothesis of homogeneity so that the expectations $\EE(E_{0,\text{opt},j})=E_{\text{mean}}$ may be assumed to be independent of $j$.

It turned out that the estimator $\sigma_{E_{0,\text{opt}}}$ has a bias depending on the correlation length of the underlying random field $E(x,y)$ in \eqref{eq:RFEnu}.
In fact, $\sigma_{E_{0,\text{opt}}}$ underestimates $\sigma_E$ for short correlation lengths, the underestimation getting more distinct for smaller correlation lengths.

As detailed in Subsection~\ref{subsec:advantages}, the effect can be explained as follows. When the correlation length gets smaller and smaller, one is entering the quasi-homo\-ge\-neous regime. The medium behaves like a deterministic one with almost constant $E$ \cite{Wu:1988}. Therefore, the estimator $\sigma_{E_{0,\text{opt}}}$ tends to zero as well. On the other hand, when the correlation length gets larger, the underestimation disappears because then the random field $R_E(x,y)$ tends to become a constant random variable $R_E$ with standard deviation $\sigma_E$, for which $\sigma_{E_{0,\text{opt}}}$ is an unbiased statistical estimator. The result is a nonlinear dependence of $\sigma_{E_{0,\text{opt}}}$ on $L_E$, as shown in Figure~\ref{fig:gammasigmae0}.

Numerical experiments showed that $\sigma_{E_{0,\text{opt}}}$ is a linear function of $\sigma_E$ with a slope depending on $L_E$, that is,
$\sigma_{E_{0,\text{opt}}} = \mathfrak{h}(L_E) \sigma_E$ for some function $\mathfrak{h}(L_E)$. It is proposed here that $\mathfrak{h}(L_E)$ can be obtained as a calibration curve through numerical experiments as follows.
First, one fixes a standard deviation $\sigma_E$ and models the dependence of the estimator \eqref{eq:sigEopt} on the correlation length, that is one tries to find a functional relationship $\sigma_{E_{0,\text{opt}}} = \mathfrak{f}(L_E)$. For this, one selects a set of correlation lengths $(L_E)_i$, $i=1,\ldots, m$. Then, for each correlation length one produces a Monte Carlo sample of signals $\boldsymbol{u}_{\text{meas}}(x_j,y_j,t)$ by means of FE-simulations, where the random field has the fixed standard deviation $\sigma_{E}$ and correlation length $(L_E)_i$. For the $i^{\text{th}}$ sample, one estimates $(\sigma_{E_{0,\text{opt}}}^2)_i$ by the procedure described above, resulting in the values
$(\sigma_{E_{0,\text{opt}}})_i = \mathfrak{f}((L_E)_i)$. The curve $\mathfrak{f}(L_E)$ is obtained through regression by means of a power function. The linear dependence of $\sigma_{E_{0,\text{opt}}}$ on $\sigma_E$ allows one to compute the slope by $\mathfrak{h}(L_E) = \mathfrak{f}(L_E)/\sigma_E$. An explicit example of the calibration procedure will be given in Section \ref{sec:numerical}.

This calibration curve can be used for any later experiment. Thus, if one has to determine the standard deviation of a new sample of the same material, one can estimate it by
\begin{equation} \label{eqn:sigEoptgleichsigE}
 {\sigma}_{E}^{*}  = \sigma_{E_{0,\text{opt}}}/\mathfrak{h}(L_E)
\end{equation}
This estimation procedure, however, is only possible if one knows $L_E$. The next subsection indicates how the correlation length $L_E$ can be estimated.

Similarly, the variance of $\nu_{0,\text{opt}}$ is estimated as:
\begin{equation}\label{eq:signuopt}
\sigma_{\nu_{0,\text{opt}}}^2 = \frac{1}{N\,n-1}\sum_{j=1}^n \sum_{k=1}^N \left (\nu_{0,\text{opt},j}^{(k)} - \nu_\text{mean}^* \right)^2.
\end{equation}
Numerical experiments showed that $\sigma_{\nu_{0,\text{opt}}}$ is a sufficiently accurate estimator of $\sigma_\nu$. In any case, in the later application the coefficient of variation of $\nu$ was taken rather small (around 1\%), so a possible bias of the
estimator for $\sigma_\nu$ was not further investigated, and
\begin{equation}\label{eqn:signuoptgleichsignu}
{\sigma}_{\nu}^{*} = \sigma_{\nu_{0,\text{opt}}}
\end{equation}
was taken as estimate for the standard deviation of Poisson's ratio.

\subsection{Estimating the correlation length}
\label{subsec:corrlength}

Since data at different sensor locations are available, one can also estimate the correlation between the sensors and from there the overall correlation length. It is assumed that the parameters in the sensor locations derive from a stationary random field obeying the autocovariance function
\begin{equation}\label{eq:corremp}
\Cov(E_{0,\text{opt},j_1},E_{0,\text{opt},j_2}) = C(r) = \sigma_{E_{0,\text{opt}}}^2 \exp\left(-\frac{r}{L_{E_0}}\right)
\end{equation}
with yet unknown correlation length $L_{E_0}$, where $r=\sqrt{(x_{j_1}-x_{j_2})^2+(y_{j_1}-y_{j_2})^2}$ is the Euclidean distance between sensors $j_1$ and $j_2$. The sensors are paired into groups with the same Euclidean distance, denoting the $l^{\text{th}}$ group of sensor pairs with Euclidean distance $r_l$ by $\nP_l$. Then the covariance at lag $r_l$ is empirically estimated by
\begin{equation*}
C_{\rm emp}(r_l) = \frac{1}{N\abs {\nP_l}-1}  \sum_{k=1}^N \sum_{[j_1,j_2]\in\nP_l} \left(E_{0,\text{opt},j_1}^{(k)} - E^*_{\text{mean}}\right)\left(E_{0,\text{opt},j_2}^{(k)} - {E}^*_{\text{mean}}\right).
\end{equation*}
The correlation length $L_{E_0}$ is then estimated by fitting the autocovariance function \eqref{eq:corremp} to the empirical covariances, resulting in the estimate
\begin{equation}\label{eq:estLE}
L_{E_0,\text{opt}}=\argmin_{L_{E_0}} \menge{ \sum_l \abs{C_{\rm emp}(r_l) - \sigma_{E_{0,\text{opt}}}^2 \exp\left(-\frac{r_l}{L_{E_0}}\right)}}.
\end{equation}
Similar to the observed bias of the estimator $\sigma_{E_{0,\text{opt}}}$, the estimator $L_{E_0,\text{opt}}$ overestimates $L_E$ when $L_E$ gets small.
Numerical experiments suggested, though, that the estimator depends on the correlation length through a functional relationship
$L_{E_0,\text{opt}} = \mathfrak{g}(L_E)$. To determine this relationship, one can use the same simulation data as used to calibrate $\mathfrak{f}(L_E)$. This time, for every $i=1,\ldots,m$, one computes $(L_{E_0,\text{opt}})_i$ and obtains the relation $(L_{E_0,\text{opt}})_i = \mathfrak{g}((L_E)_i)$. Then, by regression, one can fit the curve $\mathfrak{g}(L_E)$ and for new samples one may obtain the estimator $L_{E}^*$ by inverting $\mathfrak{g}$ and setting
\begin{equation}\label{eq:LEast}
L_{E}^* = \mathfrak{g}^{-1}(L_{E_0,\text{opt}}).
\end{equation}
Finally, $L_{E}^\ast$ is inserted in place of $L_E$ in \eqref{eqn:sigEoptgleichsigE} to get the estimate ${\sigma}_{E}^{*}$ of the standard deviation.

\section{Testing for damage}
\label{sec:damage}

This section serves to present a proposal for statistical tests for damage, indicating whether the properties of a given material have undergone some changes or not.
The starting point is a calibrated FIO model of the undamaged, random material. The response of the undamaged material is modelled through equations \eqref{eqn:linelst2d_num_sol_u}, \eqref{eqn:linelst2d_num_sol_u1_hat}, \eqref{eqn:linelst2d_num_sol_u2_hat} with wave speeds \eqref{eq:randwavespeeds} given by the random fields
\eqref{eq:RFEnu}, \eqref{eq:basicautocov}. The parameters of the random fields are calibrated by the procedure outlined in Section\;\ref{sec:Parameter_estimation}. The sensor locations remain fixed.

The strategy is as follows. Using the stochastic Fourier integral operator model of the undamaged material, a sample of signals in the sensor locations is generated. Thereby, percentiles of the distribution of certain features of the undamaged material are obtained. These data are compared with measured data of the possibly damaged material. The procedure can be cast in the form of a hypothesis test, in which the null hypothesis (undamaged material) is accepted or rejected, depending on whether the measured features remain within certain bounds or not. The generation of the Monte Carlo sample can be done in advance of the testing. Generating a large Monte Carlo sample of system responses in finitely many sensor locations through the FIO representation is computationally inexpensive.

Details will be presented in Section\;\ref{sec:numerical}. Here is a summary of the testing strategy.
\smallskip

\emph{Step 1.} A Monte Carlo sample of size $N_{\rm MC}$ of the random fields \eqref{eq:RFEnu}, \eqref{eq:basicautocov} is generated. The corresponding solutions \eqref{eqn:linelst2d_num_sol_u} are computed in the sensor locations, using the Fourier integral operator representation.
\smallskip

\emph{Step 2.} One selects features of each signal at each sensor location. This choice is critical, since this has a very strong influence on the efficiency of the test. Step 1 delivers a Monte Carlo sample of size $N_{\rm MC}$ of each feature at each sensor location.
\smallskip

\emph{Step 3.} The null hypothesis (undamaged material) is cast in the form of critical percentile ranges for the features.
\smallskip

\emph{Step 4.} Data of the specimen to be tested are acquired by measuring the signals of the sensors. Then, for each feature one computes the $p$-value with respect to the Monte Carlo sample.
\smallskip
	
\emph{Step 5.} Step 4 results in a tuple of $p$-values for each sensor location. An overall $p$-value can be aggregated by e.g. computing the overall minimum or mean. In this paper, the mean $p$-value for each sensor location was computed and the minimum of these values was taken as the overall $p$-value.
\smallskip
	
\emph{Step 6.} The null hypotheses is rejected if the overall $p$-value is smaller than a given threshold, for example, 1\% or 5\%.
\smallskip	

\emph{Remark.}
The aggregation of $p$-values is a much debated issue \cite{Heard:2018,Lehmann:2005,Vovk:2020,Wasserstein:2019}, and there are certainly more sophisticated approaches. However, the method chosen in Step 5 worked well. In fact, the performance of the test was validated in Subsection\;\ref{subsec:performance} where it was found to obey the designed acceptance and rejection rates.

\section{Numerical example}
\label{sec:numerical}

In this section, an application of the concepts to a prototypical example, a plane strain problem, will be presented. The section starts with the description of the model and the material assumptions, followed by a presentation of the proposed parameter calibration method. Next, statistical tests for damage are put to work in four scenarios. The section ends with a validation and performance analysis of the devised tests.

\subsection{The numerical model}
\label{subsec:FEM}
\begin{figure}[htb]
	\centering
	(a)\ \includegraphics[width=0.325\linewidth]{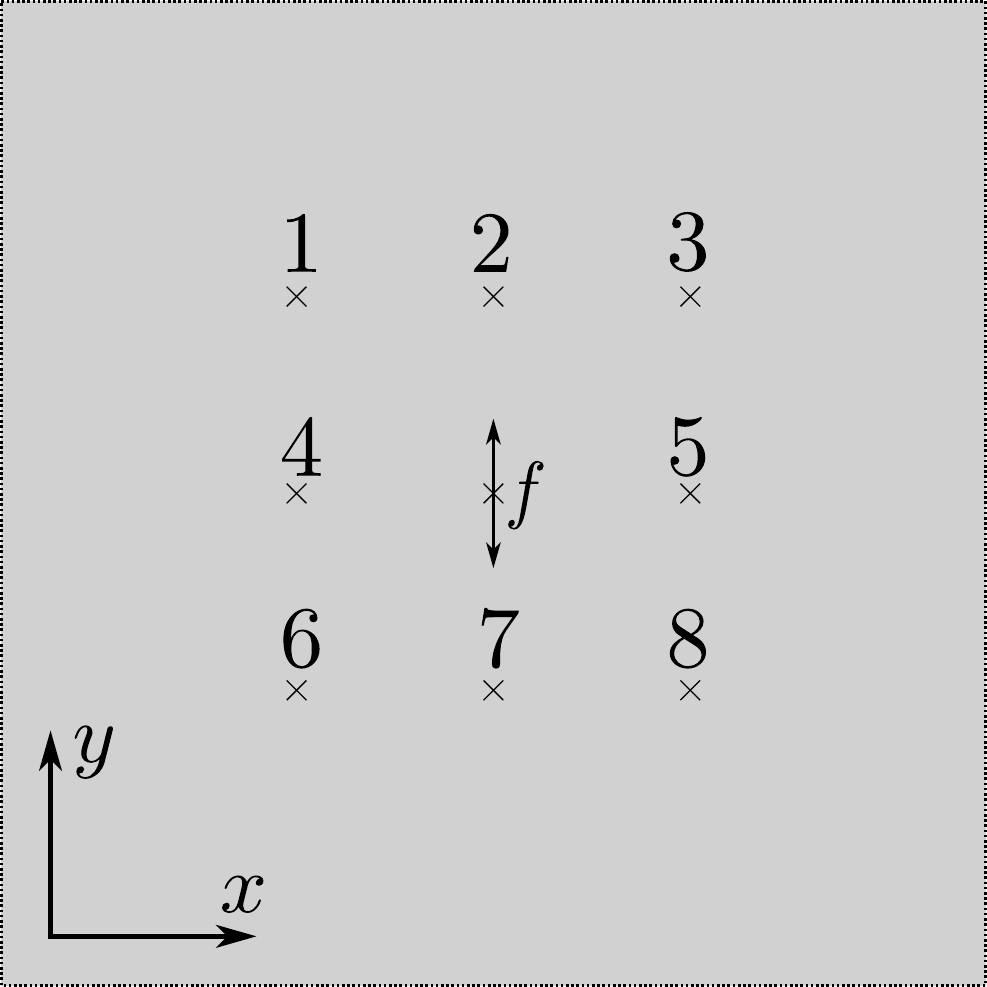}\qquad (b)\
	\raisebox{1.9cm} {\begin{tabular}{|c|c|c|}
			\hline
			$j$ & $x_j$ & $y_j$\\\hline\hline
			$1$ & $-1.17$ cm & $-1.17$ cm\\
			$2$ & $0$ cm & $-1.17$ cm\\
			$3$ & $1.17$ cm & $-1.17$ cm\\
			$4$ & $-1.17$ cm & $0$ cm\\
			$5$ & $1.17$ cm & $0$ cm\\
			$6$ & $-1.17$ cm & $1.17$ cm\\
			$7$ & $0$ cm & $1.17$ cm\\
			$8$ & $1.17$ cm & $1.17$ cm.\\\hline
	\end{tabular}}
	\caption{(a) Numerical domain with location of driving force $\boldsymbol{f}$ and position of the eight sensors; (b) coordinates of sensor positions.}
	\label{fig:paramest}
\end{figure}

The material structure is an aluminum block with the following nominal material parameters: Young's modulus is $E=70$ GPa, Poisson's ratio is $\nu=0.35$, and the density is $\rho=2.70$g/cm${}^3$. The cross section of the block was assumed to be a square of size $10\times10$ cm (henceforth referred to as the numerical domain). The geometric configuration can be seen in Figure~\ref{fig:paramest}. There are $n=8$ (virtual) sensors located around the centered force.

As in Subsection~\ref{subsec:solving_2D_model} it is assumed that the force term $\boldsymbol{f}$ is independent of $z$ and has the special form
\[
\boldsymbol{f}(x,y,z,t)= [0, g(x,y) h(t),0 ]^\sT,\quad g(x,y)=\delta(x,y), \quad h(t) =
   \left\{\begin{array}{ll}
    \sin(2\pi t), & t \geq 0,\\
    0, & t<0,
    \end{array}\right.
\]
where $\delta(x,y)$ is the Dirac delta function. The displacements $u_1(x,y), u_2(x,y)$ at point $(x,y)$ are then given by \eqref{eqn:linelst2d_num_sol_u}.

When evaluating this formula, some numerical issues have to be addressed.
Since the Fourier transform of $\delta(x,y)$ is $\widehat{\delta}(\xi,\eta)\equiv1$, the integrals in \eqref{eqn:linelst2d_num_sol_u} are not convergent.
To remedy this problem, the Dirac delta function was replaced by a regularized version, namely by
$$ \delta_\varepsilon(x,y)=\frac{1}{2\pi \varepsilon^2} \exp\big(-\frac{x^2+y^2}{2\varepsilon^2}\big),$$
where $\varepsilon$ is a small regularization parameter (chosen here as $\varepsilon = 0.1$, i.e., about three times the grid length of the numerical domain). The Fourier transform is then explicitly given by $\widehat{\delta_\varepsilon}(\xi,\eta)=\exp\big(-\varepsilon^2 (\xi^2+\eta^2)\big)$.

The numerical domain was discretized with $256\times256$ points. The considered time interval was $7\ \mu \mathrm{s}$ with $\Delta t=0.05\ \mathrm{\mu s}$. The integrals were computed by the rectangular quadrature formula for all sensor locations. For sufficient accuracy, higher order integration schemes were not required.

At $\xi=\eta=0$ one faces the additional problem that the integrands in \eqref{eqn:linelst2d_num_sol_u1_hat} and \eqref{eqn:linelst2d_num_sol_u2_hat} are not defined. This does not matter analytically since the value at one point does not change the integral. However, it is a problem for numerical evaluations. So in order to compute the solution for this point one has to go back to equation \eqref{eqn:linelst2d_full_eqn}. If one takes the Fourier transform and sets $\xi=\eta=0$ then
\begin{align*}
&\p_{tt} \vector{\widehat{u}_1(0,0,t)\\ \widehat{u}_2(0,0,t)} = \vector{0\\  h(t)},
\end{align*}
using the fact that $\widehat\delta_\varepsilon(0,0)=1$. Since the material is at rest at $t=0$, it follows that $ \widehat{u}_1({0},0,t)=0$ and $ \widehat{u}_2({0},0,t)= \int_0^t\int_0^s h(r) \d r\d s$.

The computation of the time dependent signal in all $8$ sensor locations takes approximately $4$ seconds on a workstation.

\subsubsection*{Numerical validation}

In the deterministic case of constant material parameters, the solution procedure via the Fourier integral operator formula  \eqref{eqn:linelst2d_num_sol_u} was validated by comparison with a finite element simulation of the plane strain model, solved by \eqref{eqn:linelst2d_full_eqn}, on the same domain with the same material parameters.
The same force was applied, but in non-regularized form as a point force. For the simulation, a standard Abaqus routine was used, with the implicit Euler scheme for time integration.

The (time-dependent) solution in the eight locations $(x_j,y_j)$, $j=1,\ldots, 8$ was computed at all times $t$. The evaluation of the FE-model took approximately $800$ seconds on a workstation with $5$ CPU cores used.

The comparison of the FIO-solution and the FE-solution in $x$- and $y$-direction can be found in Figures~\ref{fig:signalx0} and \ref{fig:signaly0} in the Appendix and showed satisfactory accuracy. The relative error is with respect to the maximum displacement occurring in the FE-solution in the respective directions. One should point out that there is a systematic error in the FIO-solution, since the point force was regularized. Therefore, the force term acts on an area and not at a point. That means that the FIO-signal is slightly blurred. In addition, this also means that the signal arrives at the sensor locations slightly earlier than the FE-signal. Of course, if one refines the mesh, the regularization of the point force can be more sharp and the error vanishes.

The finite element model was also used to simulate artificial data, meant to represent laboratory measurements in the next subsections.

\subsection{Parameter estimation}
\label{subsec:parest}

This subsection serves to demonstrate the calibration procedure of Section\;\ref{sec:Parameter_estimation} in the numerical example. Artificial data were generated by the FE-model of subsection\;\ref{subsec:FEM}. Young's modulus and Poisson's ratio were entered as random fields of the form \eqref{eq:RFEnu} with autocovariance functions \eqref{eq:basicautocov}. Input parameter values were $E_{\rm mean} = 70$ GPa, $\nu_{\rm mean}=0.35$, $\sigma_E=3.5$ GPa ($5\%$ of its nominal value), $\sigma_\nu=0.005$, and  $L_E=L_\nu=3$ cm. For the parameter calibration, the output of the FE-calculation was extracted at the eight sensor locations (Figure\;\ref{fig:paramest}). This output served as measured data for the FIO-based calibration procedure.

Concerning the estimation of the nominal values $E_\text{mean}$ and $\nu_\text{mean}$, three scenarios were tested, according to Table~\ref{tbl:parameterestimation}.
\begin{itemize}
\item[(I)] Deterministic case: no random fields present, Young's modulus and Poisson's ratio held constant at their nominal values $E_\text{mean}$ and $\nu_\text{mean}$.
No repetitions in formula\;\eqref{eq:EstarNustar}, $N=1$.
\item[(II)] Full random field model for the elastic parameters, no repetitions in formula\;\eqref{eq:EstarNustar}, $N=1$.
\item[(III)] Full random field model for the elastic parameters, means estimated by formula\;\eqref{eq:EstarNustar} with $N=10$.
\end{itemize}
The minimization of \eqref{eq:argmin} was done by means of the Nelder-Mead algorithm \cite{Hendrix:2010,Nelder:1965}.
To show that the estimate works well even with a bad initial guess for the algorithm, $E_{\text{ini}}=50$ GPa and $\nu_{\text{ini}}=0.3$ was chosen. The results of the optimization procedure as well as the computing time can be found in Table~\ref{tbl:parameterestimation}. The optimization algorithm converged after approximately 30 iterations.
\begin{table}[htb]	
	\centering
	\begin{tabular}{|r|c|c|c|c|}\hline
		&$E_{\text{mean}}^*$ [GPa] & $\nu_{\text{mean}}^*$ & N & time [s]\\\hline
		(I)&70.41  &  0.3490 &1 & 216 \\ \hline
		(II) &69.56  &  0.3470&1 & 214 \\ \hline	
		(III)&70.05  &  0.3473&10 & 2155 \\ \hline
	\end{tabular}
	\caption{Results of the optimization for calibrating the nominal values.}
	\label{tbl:parameterestimation}
\end{table}

\begin{figure}[htb]
	\centering
	\includegraphics[width=0.75\linewidth]{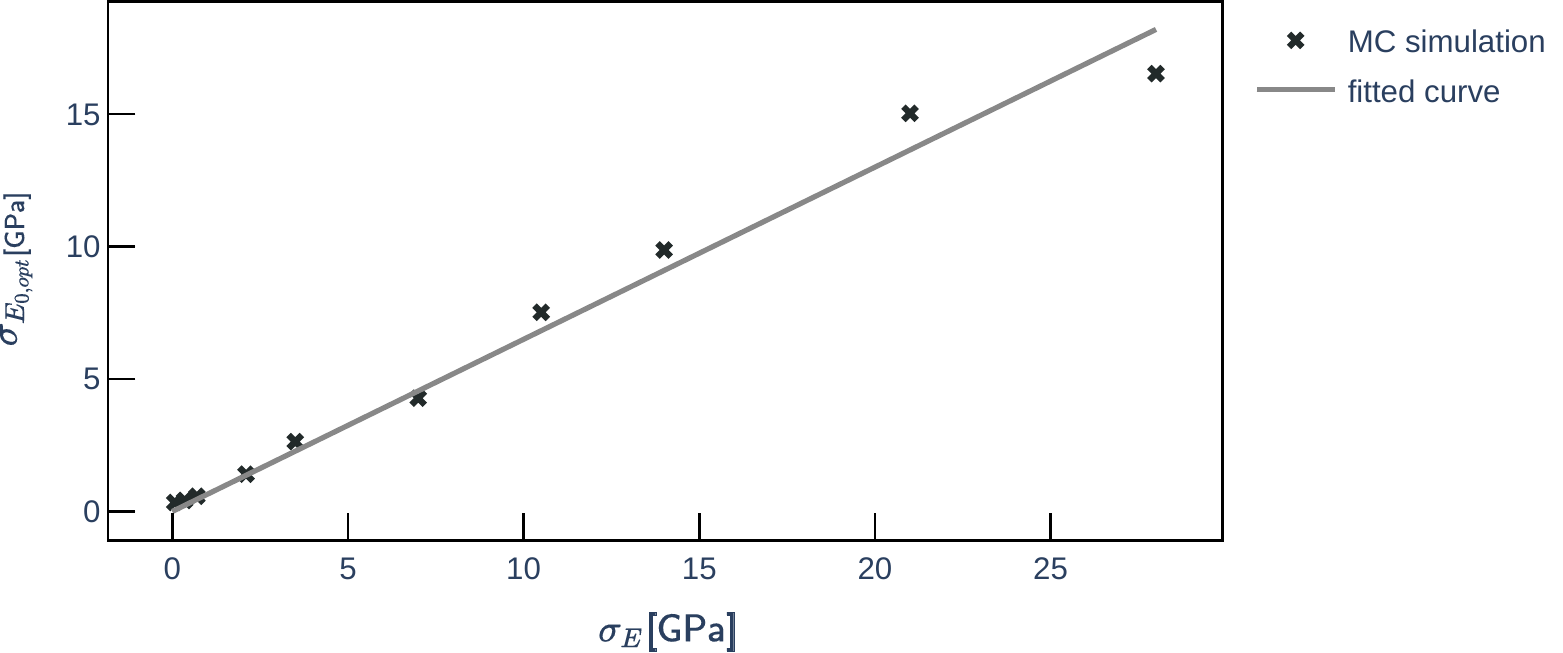}
	\caption{Linear dependence of $\sigma_{E_0,\text{opt}}$ on $\sigma_E$  with fixed correlation length $L_E=3$ cm.}
	\label{fig:sigmasigma}
\end{figure}

\begin{figure}[htb]
	\centering
	\includegraphics[width=0.75\linewidth]{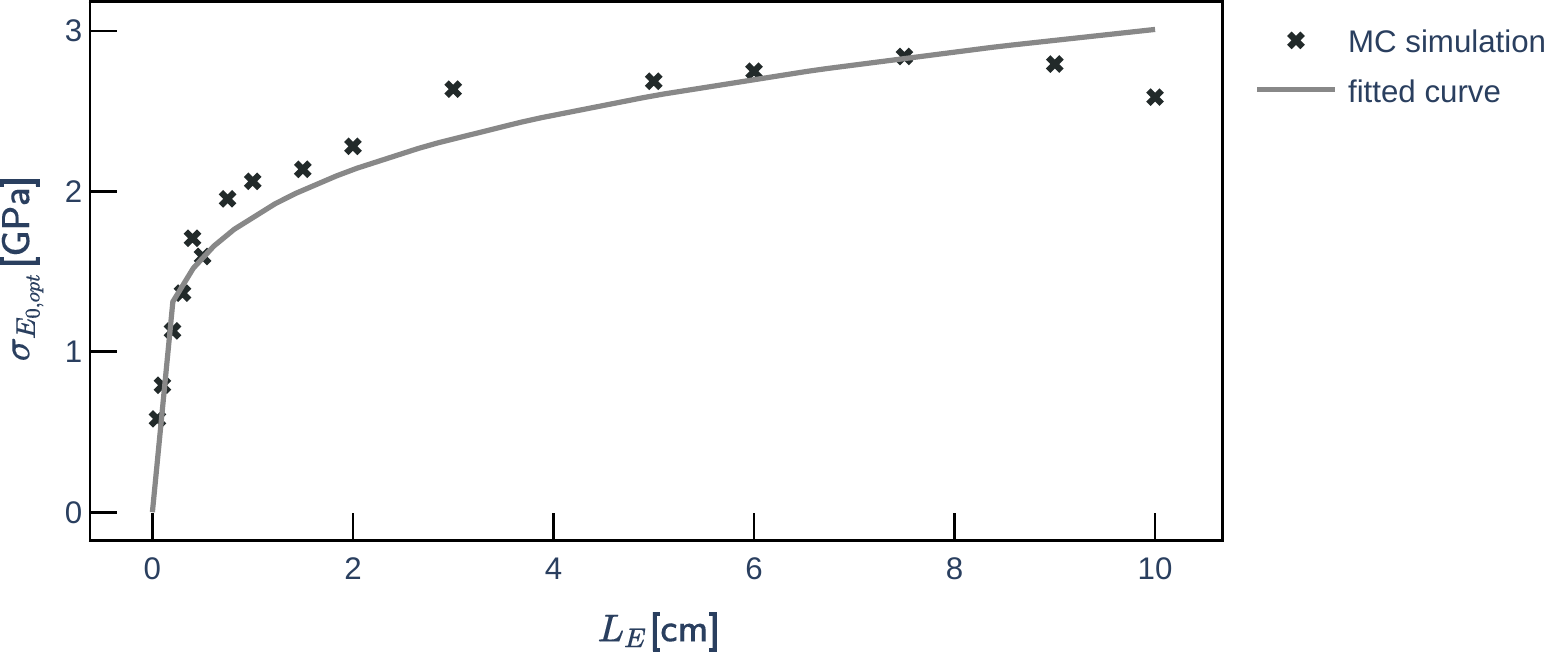}
	\caption{Dependence of $\sigma_{E_0,\text{opt}}$ on $L_E$ for fixed $\sigma_E = 3.5$ GPa. The fitted curve is given by $\mathfrak{f}(L_E) = \beta_0 L_E^{\beta_1}$ with $\beta_0 = 1.8415$ and $\beta_1 = 0.2132$.}
	\label{fig:gammasigmae0}
\end{figure}

In order to estimate the standard deviations, formulas \eqref{eq:sigEopt} and \eqref{eq:signuopt} were used with $N=100$ finite element simulations to compute $\sigma_{E_{0,\text{opt}}}$ and $\sigma_{\nu_{0,\text{opt}}}$. As indicated in \eqref{eqn:signuoptgleichsignu}, $\sigma_{\nu_{0,\text{opt}}}$ could be taken as estimator $\sigma_\nu^*$ of the standard deviation of Poisson's ratio. (The present simulation resulted in the estimate $\sigma_\nu^*= 0.0036$.) In Subsection\;\ref{subsec:stds}, the bias of $\sigma_{E_{0,\text{opt}}}$ has been pointed out. As noted there, further numerical experiments had shown that it can be described by a relation of the form $\sigma_{E_{0,\text{opt}}} = \mathfrak{h}(L_E) \sigma_E$. For example, the linear dependence of $\sigma_{E_{0,\text{opt}}}$ on $\sigma_E$ at fixed correlation length $L_E=3$ cm can be seen in Figure\;\ref{fig:sigmasigma}.
The figure was obtained by Monte Carlo simulations of sample size $N=100$ in formula \eqref{eq:sigEopt} for each $\sigma_E$ from the list
$$
\menge{0.001, 0.005, 0.01, 0.03, 0.05, 0.10, 0.15, 0.20, 0.30, 0.40}\cdot 70\mbox{\;  GPa}.
$$
The calibration curve $\mathfrak{f}(L_E)$ described in Subsection\;\ref{subsec:stds} can be similarly approximated by Monte Carlo simulation. For this purpose, formula \eqref{eq:sigEopt} was evaluated with a sample size $N=100$, for each $L_E$ from the list
$$
\menge{0.05,0.1,0.2, 0.3, 0.4, 0.5, 0.75, 1, 1.5, 2, 3, 5, 6, 7.5, 9, 10}\mbox{\ cm},
$$
at fixed $\sigma_E = 3.5$ MPa. The coefficient function was obtained in the form of a power function $\mathfrak{f}(L_E) = \beta_0 L_E^{\beta_1}$ by regression through the data. Data and fitted function are depicted in Figure\;\ref{fig:gammasigmae0}. Finally, one computes
$\mathfrak{h}(L_E) = \mathfrak{f}(L_E)/3.5$.
So far, if the correlation length is assumed to be known, one can obtain the estimate $\sigma_E^*$ of the standard deviation of Young's modulus through formula\;\eqref{eqn:sigEoptgleichsigE}.

\begin{figure}[htb]
	\centering
	\includegraphics[width=0.75\linewidth]{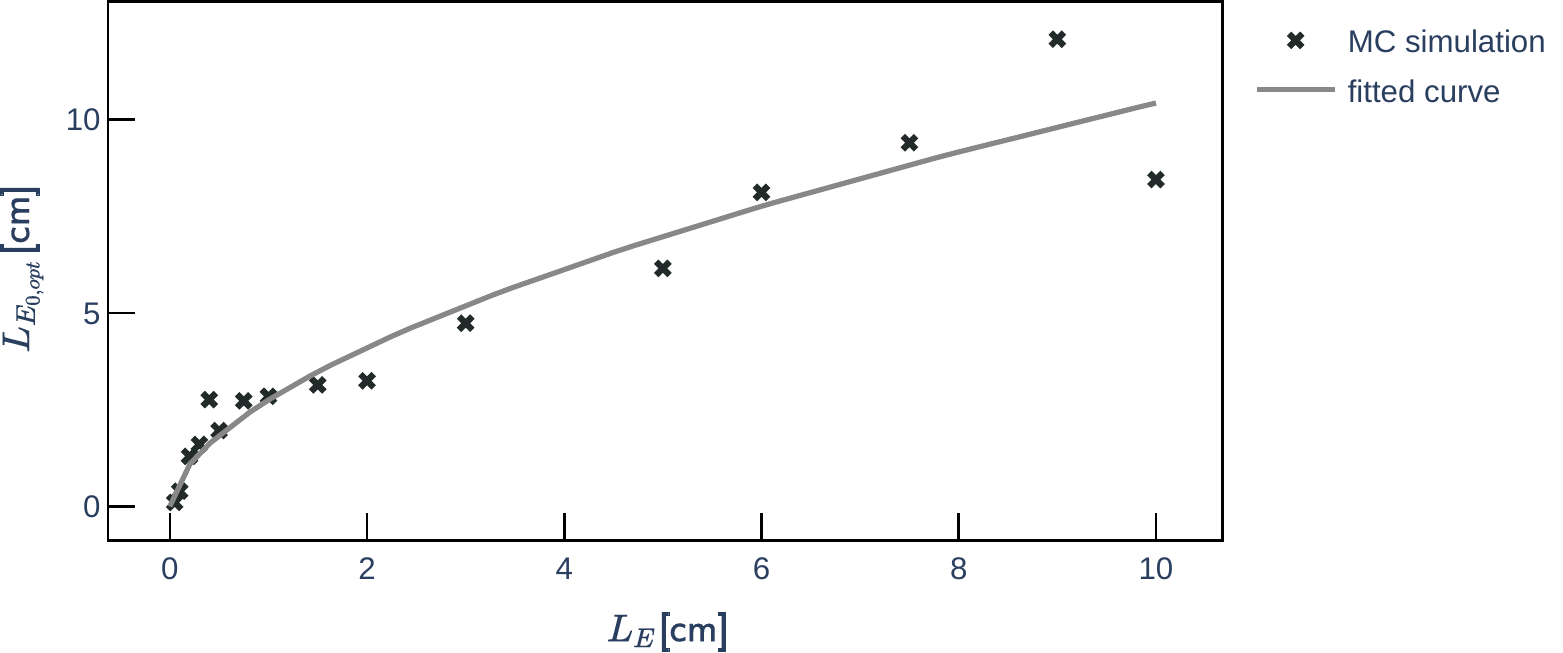}
	\caption{Dependence of $L_{E_{0,\text{opt}}}$ on $L_E$, simulation and fitted calibration curve $\mathfrak{g}(L_E) = \gamma_0 L_E^{\gamma_1}$. The obtained regression parameters were $\gamma_0 = 2.7503$ and $\gamma_1 = 0.5790$.}
	\label{fig:gammagamma}
\end{figure}

It remains to estimate the correlation lengths. The procedure has been outlined in Subsection\;\ref{subsec:corrlength}. To obtain the empirical covariances at different distances, the sensors were paired as shown in Table~\ref{table:pairs}.
\begin{table}[htb]
	\begin{center}
		\begin{tabular}{|c|c|c|c|c|c|}
			\hline
			\multicolumn{6}{|c|}{Sensor pairs ($\nP$)}                               \\ \hline
			$\nP_0$    &  $\nP_1$  &  $\nP_2$  &  $\nP_3$  &  $\nP_4$  &  $\nP_5$  \\ \hline
			$[1 , 1]$   & $[1 , 2]$ & $[2 , 4]$ & $[1 , 3]$ & $[1 , 5]$ & $[1 , 8]$ \\
			$[2 , 2]$   & $[1 , 4]$ & $[2 , 5]$ & $[1 , 6]$ & $[1 , 7]$ & $[3 , 6]$ \\
			$[3 , 3]$   & $[2 , 3]$ & $[4 , 7]$ & $[3 , 8]$ & $[2 , 6]$ &           \\
			$[4 , 4]$   & $[3 , 5]$ & $[5 , 7]$ & $[6 , 8]$ & $[2 , 8]$ &           \\
			$[5 , 5]$   & $[4 , 6]$ &           &           & $[3 , 4]$ &           \\
			$[6 , 6]$   & $[5 , 8]$ &           &           & $[3 , 7]$ &           \\
			$[7 , 7]$   & $[6 , 7]$ &           &           & $[4 , 8]$ &           \\
			$[8 , 8]$   & $[7 , 8]$ &           &           & $[5 , 6]$ &           \\ \hline\hline
			\multicolumn{6}{|c|}{Distance between sensors $(r)$}                     \\ \hline
			$r_0$     &   $r_1$   &   $r_2$   &   $r_3$   &   $r_4$   &   $r_5$   \\ \hline
			\ \ 0 cm \ \ &  1.17 cm  &  1.65 cm  &  2.34 cm  &  2.62 cm  &  3.31 cm  \\ \hline
		\end{tabular}
		\vspace*{-1.5em}
	\end{center}
	\caption{Pairing of sensors and corresponding distances.}
	\label{table:pairs}
\end{table}
Fitting the autocovariance function \eqref{eq:corremp} results in the estimator \eqref{eq:estLE}, which was observed to be biased in Subsection\;\ref{subsec:corrlength}, obeying a functional relation
$$
L_{E_{0,\text{opt}}} = \mathfrak{g}(L_E).
$$
To obtain this relation, the procedure of fitting the correlation length was repeated with sample size $N = 100$ for each $L_E$ from
$$
\menge{0.05,0.1,0.2, 0.3, 0.4, 0.5, 0.75, 1, 1.5, 2, 3, 5, 6, 7.5, 9, 10}\mbox{\ cm}.
$$
Again, a power function of the form $\mathfrak{g}(L_E) = \gamma_0 L_E^{\gamma_1}$ showed a satisfactory fit. The data and fitted function can be seen in Figure\;\ref{fig:gammagamma}.

Finally, an estimate $L_E^*$ of the correlation length $L_E$ can be achieved through formula \eqref{eq:LEast}. The results of the present numerical simulation are summarized in Table\;\ref{tbl:estimationSTDandCORR}.
\begin{table}[htb]	
	\centering
	\begin{tabular}{|c|c|c|c|c|c|}\hline
		$L_E$ [cm] & $L_{E_{0,\text{opt}}}$ [cm] & $L_E^*$ [cm] & $\sigma_E$ [GPa] & $\sigma_{E_{0,\text{opt}}}$ [GPa] & $\sigma_E^\ast$ [GPa]\\\hline
		3.0  &  4.7385 & 2.5589 & 3.5 & 2.3275 & 3.6207 \\ \hline
	\end{tabular}
	\caption{Results of the optimization for calibrating the nominal values.}
	\label{tbl:estimationSTDandCORR}
\end{table}

It can be observed that the obtained estimates of the correlation length and the standard deviation are satisfactory, but not as accurate as the estimates of the nominal values in Table\;\ref{tbl:parameterestimation}. Larger sample sizes in the determination of the calibration curves could improve the estimates, but at the price of considerably increased computational cost. The tests for damage in Subsection\;\ref{subsec:damagedet} are based on the nominal value $E_{\rm mean}$ of Young's modulus which can be accurately estimated even by a single measurement, with $N = 1$ in formula\;\eqref{eq:EstarNustar}.

\subsection{Damage detection}
\label{subsec:damagedet}

As outlined in Section\;\ref{sec:damage}, the testing strategy is based on a stochastic Fourier integral operator model of the undamaged material, which is used to generate a sample of signals in the sensor locations. This sample serves as a realization of the null hypothesis. Measured data of the possibly damaged material are tested against this sample. As before, the parameter values for the undamaged state were assumed to be $E_\text{mean}=70$ GPa, $\nu_\text{mean}=0.35$, $\sigma_E = 3.5$ GPa, $\sigma_{\nu}=0.005$, and $L_E = L_\nu = 3$ cm. These data were used to generate the random fields \eqref{eq:RFEnu}, \eqref{eq:basicautocov} producing the wave speeds \eqref{eq:randwavespeeds}, which entered the Fourier integral operator representation \eqref{eqn:linelst2d_num_sol_u} of the displacements. The coefficients of variation were so small that no problems with the square roots in \eqref{eq:randwavespeeds} could arise, and further cut-offs were not needed.

In practice, estimated values would be entered in \eqref{eq:RFEnu}, \eqref{eq:basicautocov} according to Section\;\ref{sec:Parameter_estimation}. The purpose of this subsection is to demonstrate how the test procedure works, so this step was skipped here.

The finite element model was used to generate artificial responses of the material having undergone various kinds of damage, represented through changes in the nominal value $E_\text{mean}$, or through geometric deterioration. The standard deviations and correlation lengths were kept fixed throughout.

The null hypothesis was that the material is in the undamaged state, i.e., it has the parameters described above. As initial step, a sample of size $N = 10000$ of system responses was generated, distributed according to the null hypothesis, by means of the FIO-solution
operator with wave speeds \eqref{eq:randwavespeeds}.

In order to test the material one has to choose characteristic features of the signal. This choice is quite critical, since a bad choice of features can lead to a poor distinction between damaged and undamaged material.

Note that for the tests, a full record of the measured signal (duration $7\ \mu \mathrm{s}$) is available. Using the information captured in the Fourier spectrum of the signal suggests itself, because the amplitude of a selected frequency relates to attenuation, while the phase is directly related to propagation speed. Experiments with various sets of frequencies of the Fourier spectrum were undertaken.
For the present purpose, the amplitude and the phase angle of the two most dominant frequencies of the recorded signal turned out to work best.

As described in
Subsection\;\ref{subsec:FEM}, there were eight sensor locations with signals in $x$- and $y$-direction. The four signals in $x$-direction at sensors 2, 4, 5, 7 were close to zero. (Due to the special excitation, no pressure waves arrive at sensors 4 and 5, and no shear waves arrive at sensors 2 and 7.) Thus 12 signals remained for the analysis. Of each signal $2\times2$ features were extracted:  the phase angle and amplitude at the first two nonzero frequencies in the DFT-spectrum (here $\omega_1 = 2\pi/7 \approx 0.9$ MHz and $\omega_2 = 4\pi/7 \approx 1.8$ MHz). Thus a total number of 48 features were used for damage detection.

Three tests with the following decision procedures were implemented:
\begin{enumerate}
	\item [(I)]  A left-sided test for the phase angle of the signal, detecting small phase angles (corresponding to late arrival of the signal). For this test, the phase angles need to be ordered: The phase angle space was determined by the mean phase angle of the FIO sample plus/minus $180$ degrees.
	
	For each sensor location, the $p$-value for each single feature (i.e. the phase angle of one signal for one frequency) was determined, and the average $p$-value (over all features) at each sensor location was calculated. The overall $p$-value was taken as the minimum over all sensor locations.
	\item [(II)] A right-sided test  for the phase angle of the signal, detecting large phase angles, using the same principle as in (I).
	\item [(III)] A two-sided test for the amplitude of the dominant frequencies. The aggregated $p$-value (over features and sensor locations) was computed as in (I).
\end{enumerate}

These tests were applied to the following scenarios:

\begin{itemize}
	\item[S1:] The material is having the desired properties (undamaged state).
	
	\item[S2:]  The material is not having the desired properties: $E_\text{mean}$ is too small (60 MPa).
	
	\item[S3:]  The material is not having the desired properties: $E_\text{mean}$ is too large (80 MPa).
	
	\item[S4:]  The material is having the desired properties, but is suffering a crack. The crack is modelled by a small region in the FE-model having a very small Young's modulus. The modelled crack lies between the applied point force and Sensor 3.
\end{itemize}
\begin{figure}[htb]
	\centering
	(a)\includegraphics[width=0.45\linewidth]{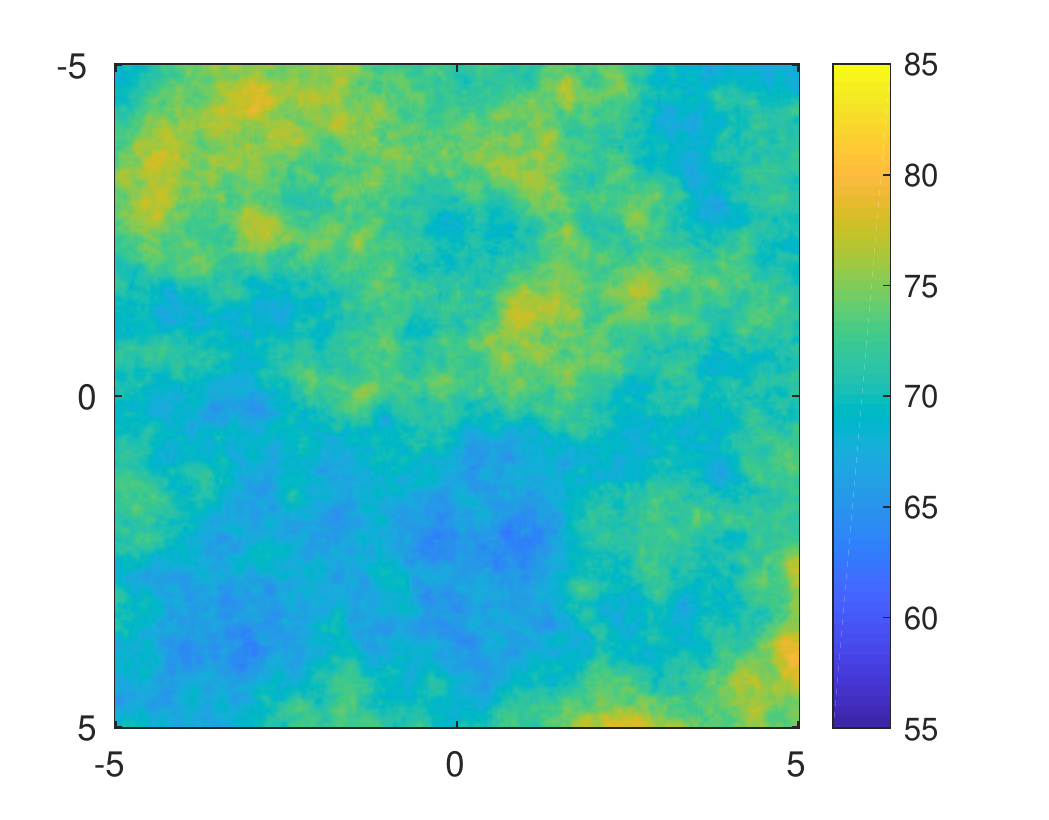}
	(b)\includegraphics[width=0.45\linewidth]{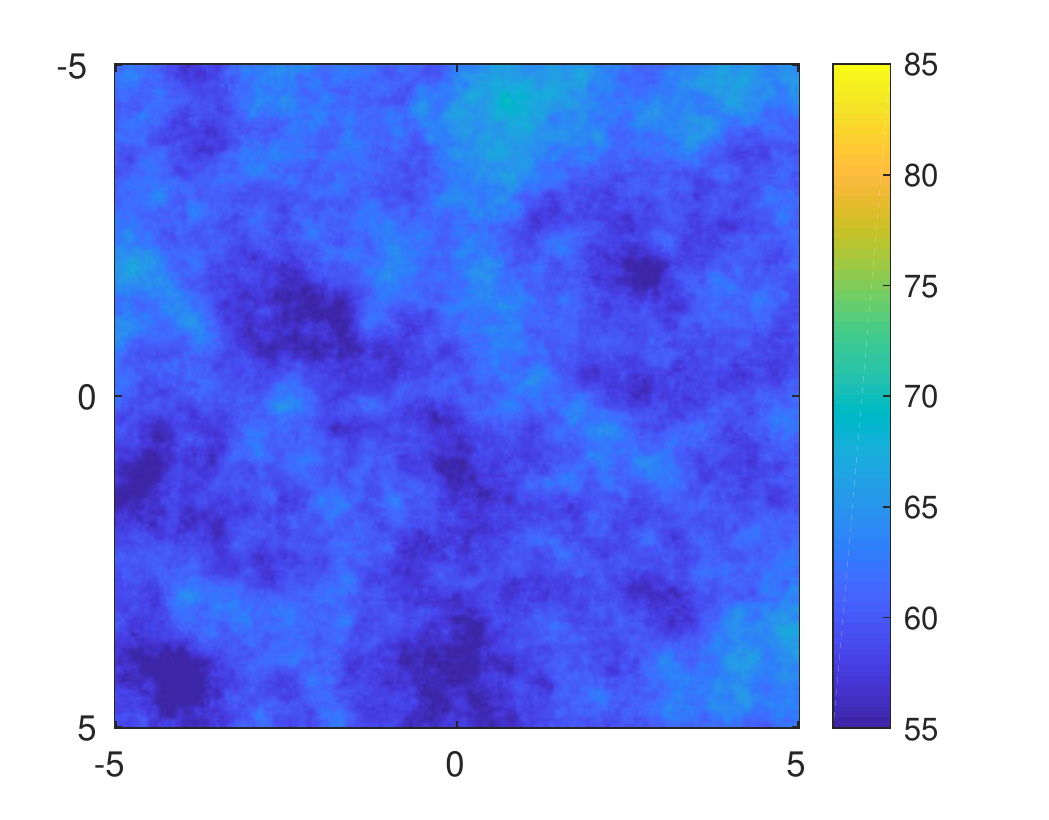}
	(c)\includegraphics[width=0.45\linewidth]{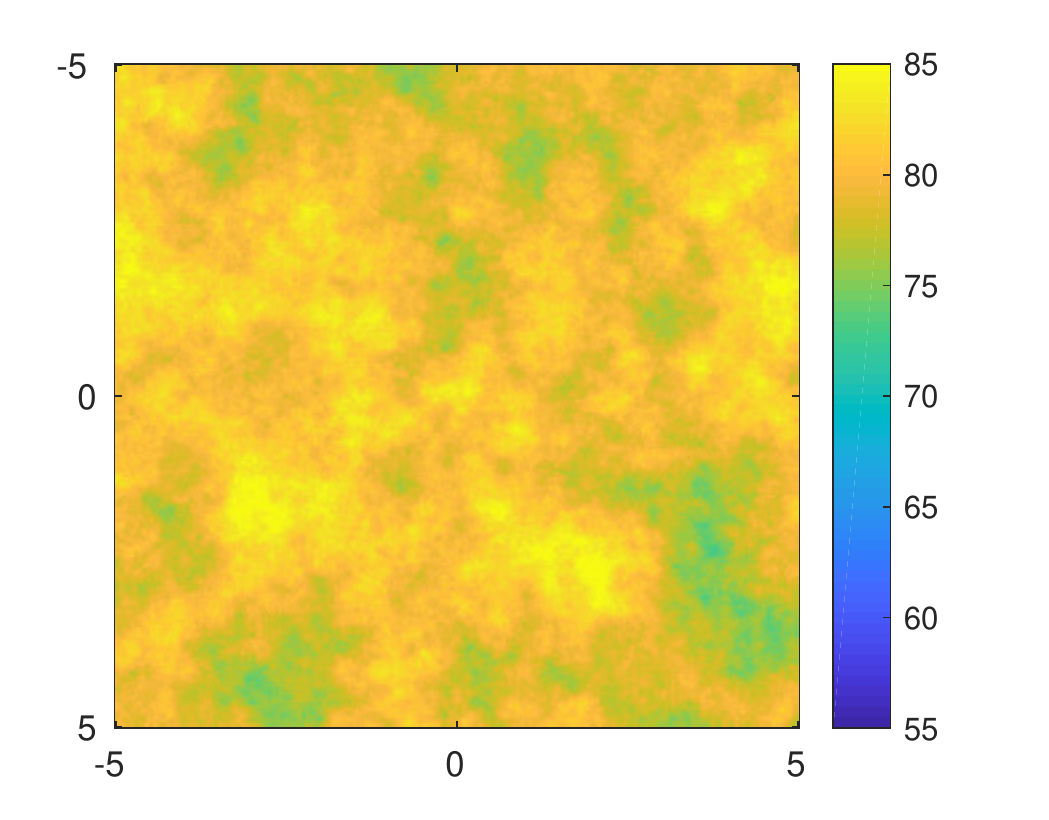}
	(d)\includegraphics[width=0.45\linewidth]{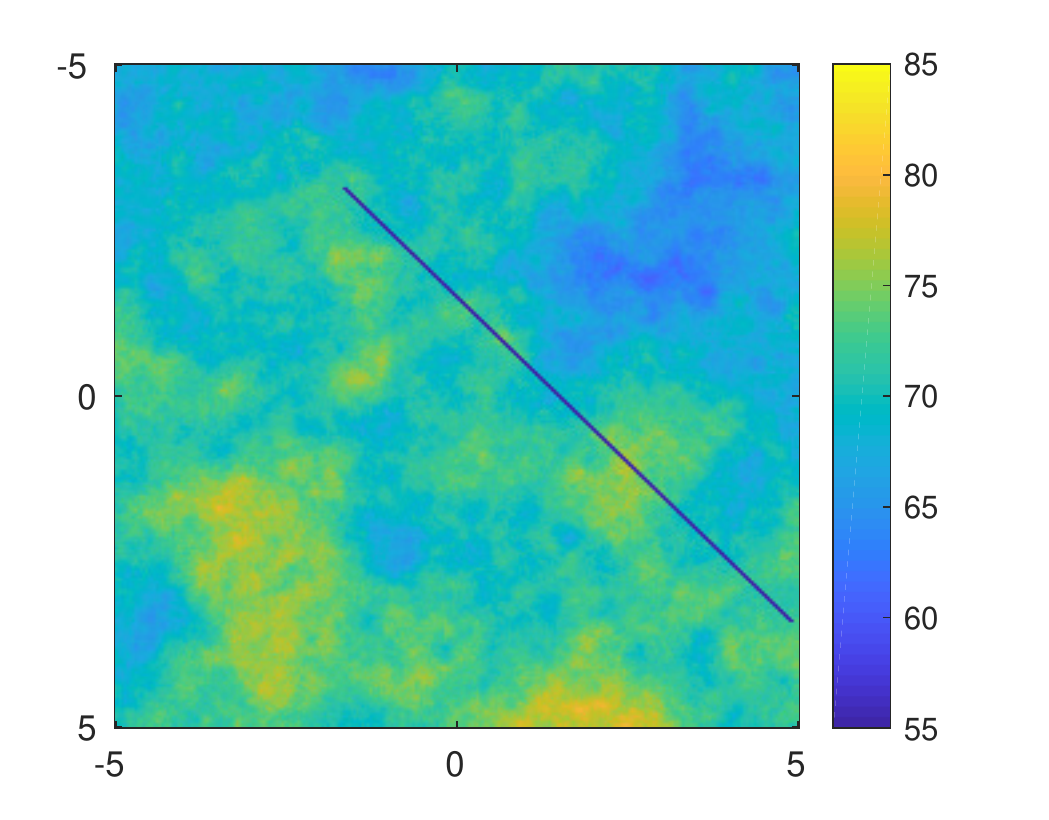}
	\caption{Realizations of the random field $E(x,y)$ (Young's modulus) in Scenarios S1--S4. (a) Scenario S1 (undamaged); (b) Scenario S2 ($E_\text{mean}$ is too small); (c) Scenario S3 ($E_\text{mean}$ is too large); (d) Scenario S4 (crack). The color scale is in GPa.}
	\label{fig:cases}
\end{figure}
Figure\;\ref{fig:cases} shows realizations of the random field $E(x,y)$ in Scenarios S1--S4. A confidence level of 99\% was adopted. This meant rejection of the null hypothesis for $p$-values below $\alpha = 1\%$. (In Subsection\;\ref{subsec:performance}, $\alpha = 5\%$ will be considered for comparison as well.) The confidence bounds of the features corresponding to the chosen $\alpha$-levels were determined through the sample of size $10000$ of simulated stochastic FIO-solutions of the undamaged material. For a typical random field realization the tests showed the following results.
\subsubsection*{Test (I)}
\begin{enumerate}
		\item [S1:]  The computed $p$-value lies well above the chosen $\alpha$-level, which means that undamaged material will be accepted as such.
		\item [S2:]  The $p$-value is well below 1\%, the phase angles are to the left of the lower confidence bounds. The test detects that $E_{\rm mean}$ is small. In fact, for small $E_{\rm mean}$, the wave speed is low, the signal arrives later than in the undamaged state, and this corresponds to small phase angles.
		\item [S3:]  The $p$-value is close to one. The left-sided test does not reject the null hypothesis. In fact, Young's modulus is high, the wave speed is large, the signal arrives early and the phase angles are actually large.
		\item [S4:] The crack has no significant influence on the phase angles.
	\end{enumerate}
\subsubsection*{Test (II)}
	This shows the same result as Test (I) with the exception that Scenario S2 has a large $p$-value and Scenario S3 has a low $p$-value. This is due to the fact that this is a right-sided test.
\subsubsection*{Test (III)}
	\begin{enumerate}
		\item [S1:] The $p$-value of the amplitude lies above the chosen $\alpha$-level, and therefore the null hypothesis is not rejected.
		\item [S2:] The amplitude also changes with Young's modulus, but the test is not suitable for detecting this scenario.
		\item [S3:] The amplitude also changes with the Young's modulus, but the test is not suitable for detecting this scenario.
		\item [S4:] The crack influences the signal in Sensor $3$, since it is in between Sensor $3$ and the wave source. As a consequence, the amplitude changes as well. The null hypothesis is rejected in Sensor 3. In the other sensors, the behavior is varying, but the null hypothesis is usually accepted in Sensor 6, where the signal is not influenced by the crack.
	\end{enumerate}

It is interesting to note that the different behavior of Test (III) in the different sensors allows one to draw conclusions about the location of the crack. The three tests are most expressive when considered together. It also should be noted that, when repeating the tests, errors of the first kind (rejection of true null hypothesis) and of the second kind (acceptance of false alternative) do occur. However, as will be seen in the next subsection, the error rate of the first kind is around or below the chosen $\alpha$-level.

Figures~\ref{fig:signal:fine} to \ref{fig:signal:crack} show the comparison of the Monte Carlo sample of size $10000$ (representing the null hypothesis) with the FE-computed signal, based on one typical realization of the random field of the corresponding scenario in a certain sensor location.
The histograms depict the marginal distributions of the respective features generated by the Monte Carlo sample. Also shown are the two-sided $99\%$ bounds of the amplitude and the one-sided $99\%$ upper and lower bounds of the phase angle for the indicated frequency, as well as the value of the feature obtained from the FE-computation.
The top figures always depict the time dependent displacement in the indicated direction at the indicated sensor. The gray band visualizes the collection of the solution curves from the Monte Carlo sample whose features lie in the $99\%$ confidence intervals.

The asymmetries in the histograms of the amplitudes in Figures~\ref{fig:signal:smallE} to \ref{fig:signal:crack} are caused by a strong nonlinear correlation of amplitude and phase angle at frequency $\omega_1$, resulting in an inversely U-shaped scatterplot (not shown here).

\subsection{Performance evaluation of the tests}
\label{subsec:performance}

In order to check whether the tests obey the designed rates of false negative and false positive classifications, tests were repeatedly applied to FE-simulated data with different realization of the random fields. More precisely, the tests were applied to a sample of $100$ FE-simulations of each scenario. The numbers in Table~\ref{tbl:pvalues} suggest that the tests perform as designed (at $\alpha$-levels of 5\% and 1\%, respectively).
\begin{table}[htb]\centering
	\begin{tabular}{|c|c|c|c|c|}
		\hline & Scenario 1 & Scenario 2 & Scenario 3 & Scenario 4 \\ \hline
		(I), $\alpha=5\%$  &     0      &     99     &     0      &     6      \\ \hline
		(I), $\alpha=1\%$  &     0      &     90     &     0      &     0      \\ \hline
		(II), $\alpha=5\%$  &     5      &     0      &    100     &     22     \\ \hline
		(II), $\alpha=1\%$  &     1      &     0      &     99     &     3      \\ \hline
		(III), $\alpha=5\%$ &     5      &     74     &    100     &    100     \\ \hline
		(III), $\alpha=1\%$ &     0      &     45     &     89     &     99     \\ \hline
	\end{tabular}
	\caption{Number of $p$-values below the $\alpha$-level among 100 repetitions.}
	\label{tbl:pvalues}
\end{table}
The performance of the tests, based on the $100$ FE-simulations, is visualized in Figure~\ref{fig:tests_scens} by means of a scatter plot\footnote{The points are arranged according to their $p$-values.} and a violin plot\footnote{The shaded area indicates the smoothed kernel density estimator of the $p$-values. The visualization program was taken from {MathWorks, Violin Plots for plotting multiple distributions. https://www.mathworks.com/matlabcentral/fileexchange/23661-violin-plots-for-plotting-multiple-distributions-distributionplot-m}, 2017 (accessed 12 July 2018).}.

As a conclusion, one can say that the performance analysis shows that the proposed procedure is capable of detecting damaged material, and furthermore the three tests are capable of distinguishing between the scenarios.

\subsubsection*{The influence of $L_E$ and $\sigma_E$}

As was seen from the parameter estimation, the correlation length of $E$ plays a critical role in estimating the parameters. So the question arises if it also influences the $p$-values of the tests. To check this point, the tests from above were applied to an undamaged material (i.e. $E_{\text{mean}}=70$ GPa and $\sigma_{E}=3.5$ GPa) with different correlation lengths $L_E$. For each correlation length, a sample of finite element simulations with sample size $100$ was tested against the Monte Carlo sample of size $10000$ from Subsection\;\ref{subsec:damagedet}. The results can be seen in Figure~\ref{fig:tests_Ls}. One may observe that the correlation length does play a role, but does not influence the results too critically. The number of rejections can be found in Table~\ref{tbl:pvalues_of_L}.

In contrast to the correlation length, the standard deviation severely influences the measured signals. The tests from above were again applied to an undamaged material (i.e. $E_{\text{mean}}=70$ GPa and $L_{E}=3$ cm) with different $\sigma_E$. For each standard deviation value, a sample of finite element simulations with sample size $100$ was tested against the Monte Carlo sample from Subsection\;\ref{subsec:damagedet}. The results can be seen in Figure~\ref{fig:test_sigma}. Up to a coefficient of variation of $E$ of 5\%
all $p$-values were above the chosen $\alpha$-level, at least in Test (I). If the coefficient of variation was larger than $15\%$, the null hypothesis was falsely rejected for almost the whole sample of size $100$. However, since a large standard deviation can be interpreted as damaged material, too, this can be seen as a desirable property of the tests. The number of rejections can be found in Table~\ref{tbl:pvalues_of_sigma}.

\subsection{Advantages and limitations}
\label{subsec:advantages}

\subsubsection*{Advantages of the method}

The advantage of the described method is its low computational cost. The evaluation of $10000$ FIOs took approximately $11$ hours ($=10000\cdot4$ sec). If one wanted to compute a Monte Carlo sample of responses of the undamaged material of the same size by the finite element method, one could expect a computation time of $93$ days ($10000 \cdot 800$ sec).

The computation of the $p$-values is so fast that it can be done online. The computation of the features of the Monte Carlo sample of size 10000 representing the null hypothesis can be done in advance. So if one gets a measured signal one just has to compute its features and compare them with the Monte Carlo sample. This can be done in real time.

\subsubsection*{Limitations and range of applicability}

The following limitation of the presented procedure should be pointed out.
If the stochastic FIO representation has the form as described in Section\;\ref{sec:FIO_solution_3d_linelst}, it does not incorporate any internal reflections (scattering, mode conversion) of the waves due to local changes of material parameters, as was already noted in the introduction and in Subsection\;\ref{subsec:stds}. The experiments subsumed in Table~\ref{tbl:pvalues_of_sigma} suggest that the applicability of the method is mostly influenced by the coefficient of variation of the material parameters, which should be sufficiently small.

In order to give a rough quantification, also with respect to frequency, wavelength, correlation length and propagation distance, one may invoke the classification of scaling regimes in \cite[Section 5.1]{Fouque2007}\footnote{With some caution, because the quoted work refers to waves in randomly layered media.}; see also \cite{Wu:1988}. Denoting the typical frequency by $f_0$, the important parameters for the occurrence or non-occurrence of scattering are: the angular frequency $\omega_0 = 2\pi f_0$; the wave number $k_0 = \omega_0/c_0$, where $c_0$ is the typical propagation speed, the coefficient of variation $\sigma$ of the underlying random fields, the scale of heterogeneities $\ell$ (here taken as the correlation length), and the typical propagation distance $L_0$. For the present purpose, the following regimes listed in \cite{Fouque2007} are of interest:
\smallskip

\emph{The effective medium:} $k_0\ell \ll 1$, $k_0L_0\sim 1$, $\sigma \ll 1$ or $\sigma \sim 1$. The medium is quasi-homogeneous, random scattering is weak, and no backscattering occurs.
\smallskip

\emph{The weakly homogeneous regime:} $k_0\ell \sim 1$, $k_0L_0\gg 1$, $\sigma \ll 1$. The coupling between the wave and the medium is weak due to small $\sigma$. Significant scattering can only occur at large propagation distances.
\smallskip

\emph{Weakly homogeneous regime, subcase} $\sigma^2k_0L_0 \sim 1$: In this case, mode coupling and backscattering are of order one.
\smallskip

In the case of the presented numerical example, the excitation frequency was $f_0 = 1\ \mu \mathrm{s}^{-1}$, $c_0 \approx \sqrt{E/\rho} \approx 0.5\ \mathrm{cm}/\mu \mathrm{s}$ (see \eqref{eq:wavespeeds}); the precise values were $c_l = 0.65\ \mathrm{cm}/\mu \mathrm{s}$, $c_s = 0.31\ \mathrm{cm}/\mu \mathrm{s}$), $\sigma \approx 0.05$ (equal to 0.05 for $E$ and less than 0.0015 for $\nu$), $\ell = 3$ cm, $L_0\approx 10$ cm. For these values, one obtains $k_0\approx 10\ \mathrm{cm}^{-1}$ and thus
\[
   k_0\ell \approx 30,\quad k_0L_0 \approx 100,\quad \sigma \approx 0.05.
\]
This puts the setup into the weakly homogeneous regime; at the considered propagation distance, at most weak scattering is expected.

With regard to the applicability of the hypothesis tests, Table~\ref{tbl:pvalues_of_L} shows the results of experiments with $\ell= L_E$ between 0.05 cm and 20 cm. The smaller $\ell$, the smaller $k_0\ell$ and the closer one gets to the quasi-homogeneous regime. This confirms the observation that in the considered frequency range, the correlation length has little influence on the hypothesis tests.

On the contrary, the coefficient of variation $\sigma$ does have an influence. In fact, for $\sigma \approx 0.1$, the subcase $\sigma^2k_0L_0 \sim 1$ is reached, where scattering of order one is to be expected. This in turn is confirmed by Table~\ref{tbl:pvalues_of_sigma}, which shows failure of the hypothesis tests beyond this value.

Summarizing, it is expected that the FIO approximation is of acceptable accuracy as long as one remains in the quasi-homogeneous or the weakly homogeneous regime and the coefficient of variation is significantly smaller than one. This holds unrestrictedly for the hypothesis tests. For parameter estimation, calibration curves as described in Section~\ref{sec:Parameter_estimation} and Subsection~\ref{subsec:parest} are required. Finally, the dimensions of the domain beyond the sensors must be chosen large enough so that the signal does not reach the domain's boundary in the measuring period (see the remark at the end of Subsection~\ref{subsec:solving_2D_model}).

\section{Conclusion}

This article presented a Fourier integral operators based approach to modelling wave propagation in random linearly elastic materials. This FIO representation of the solution to the equations of motion is numerically fast, if one is only interested in the time-dependent solution in certain locations. For example, this is the case if one needs to compare the solution with sensor data in certain locations. This efficient simulation procedure can be used for identification of the material parameters. Even if the material parameters are randomly perturbed by a Gaussian random field, the estimation of the nominal value is stable and accurate. Furthermore, using a calibration curve one can also estimate the stochastic parameters of the random fields, such as standard deviation or correlation length.

The obtained material parameters can be used for constructing a stochastic FIO, describing wave propagation in the material with its given random properties. Since the evaluation of the FIO is fast, large Monte Carlo samples of the stochastic FIO solution can be generated with relatively small effort. This sample can be used to design hypothesis tests, to determine if the material has the desired properties or not (or if the properties have changed over time). In a numerical example, four scenarios were considered: the material has the assumed properties; the material is stiff; the material is soft; a crack is present. For the decision procedures of the designed tests certain features of the time-dependent solution were considered: the phase angles and the amplitudes of the most dominant frequency of the signals. It was shown that the three implemented tests are capable of distinguishing between the four scenarios. Furthermore, it was also shown by additional simulations that the implemented tests behave as designed with regard to the error rates of false classifications.

\section*{Acknowledgements}
The authors gratefully acknowledge support by The Austrian Science Fund (FWF), grant P-27570-N26. Thanks are due to the valuable comments of both referees. Especially the detailed criticism of one of the referees led to inclusion of Subsection \ref{subsec:advantages}
and to an overall improvement of the manuscript.

\begin{appendix}
\section*{Appendix: Figures and Tables}

\begin{landscape}
		\begin{figure}[htb]
			\centering
			\includegraphics[width=0.91\linewidth]{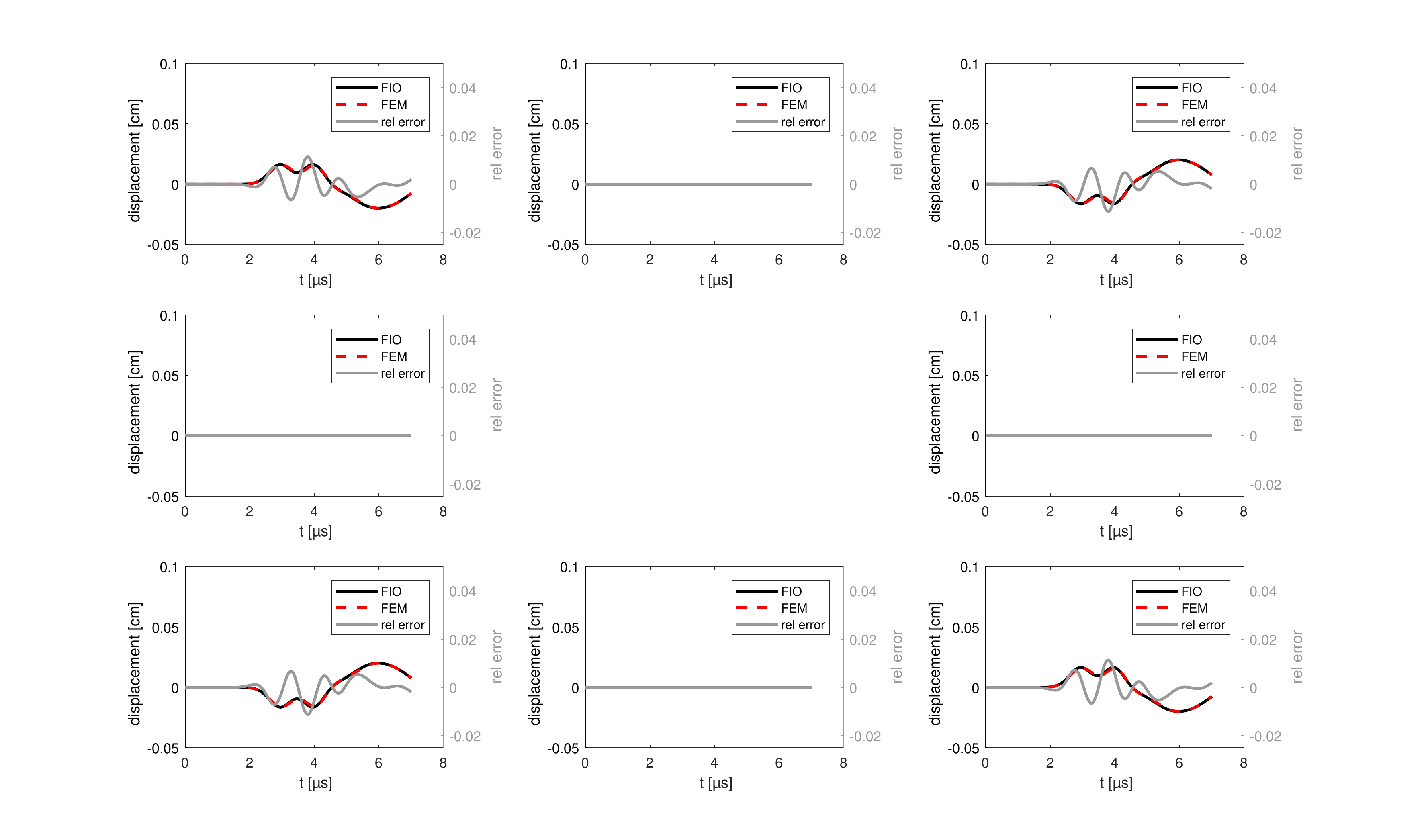}
			\caption{Comparison of FE-solution with FIO-solution, computed sensor signals -- displacement in $x$-direction. The windows are placed according to the sensor positions depicted in Figure\;\ref{fig:paramest} (top row: sensors 1--3; middle row: sensors 4 and 5; bottom row: sensors 6--7). The gray line is the relative error between FE and FIO solution with respect to the maximum displacement in the FE-solution in $x$-direction. Note that signal and error are plotted in different scales.}
			\label{fig:signalx0}
		\end{figure}
\end{landscape}

\begin{landscape}
		\begin{figure}[htb]\centering
			\includegraphics[width=0.91\linewidth]{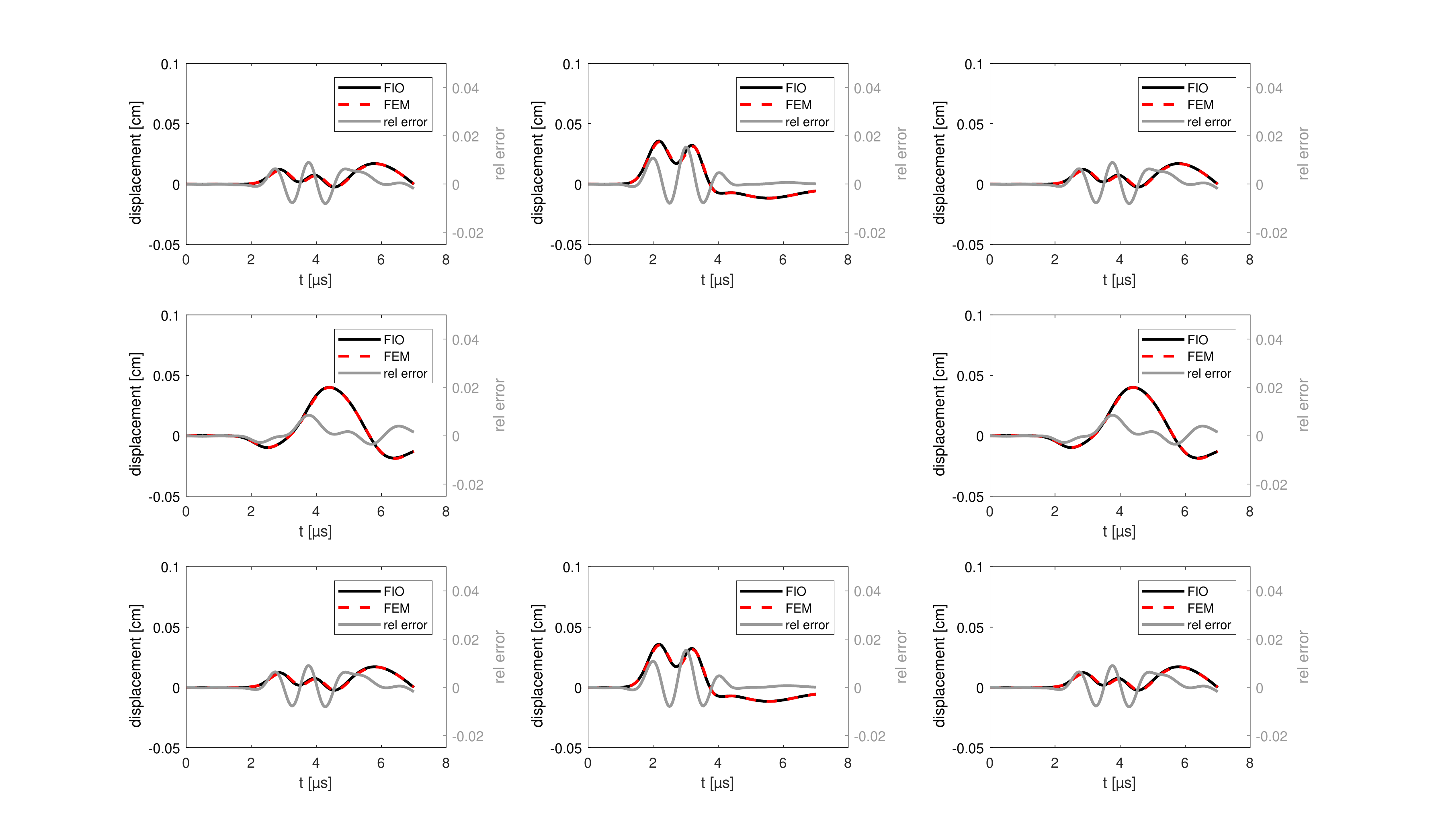}
			\caption{Comparison of FE-solution with FIO-solution, computed sensor signals -- displacement in $y$-direction. The windows are placed according to the sensor positions depicted in Figure\;\ref{fig:paramest} (top row: sensors 1--3; middle row: sensors 4 and 5; bottom row: sensors 6--7). The gray line is the relative error between FE and FIO solution with respect to the maximum displacement in the FE-solution in $y$-direction. Note that signal and error are plotted in different scales.}
			\label{fig:signaly0}
		\end{figure}
\end{landscape}

\begin{figure}[htb]\centering
		\includegraphics[width=\linewidth]{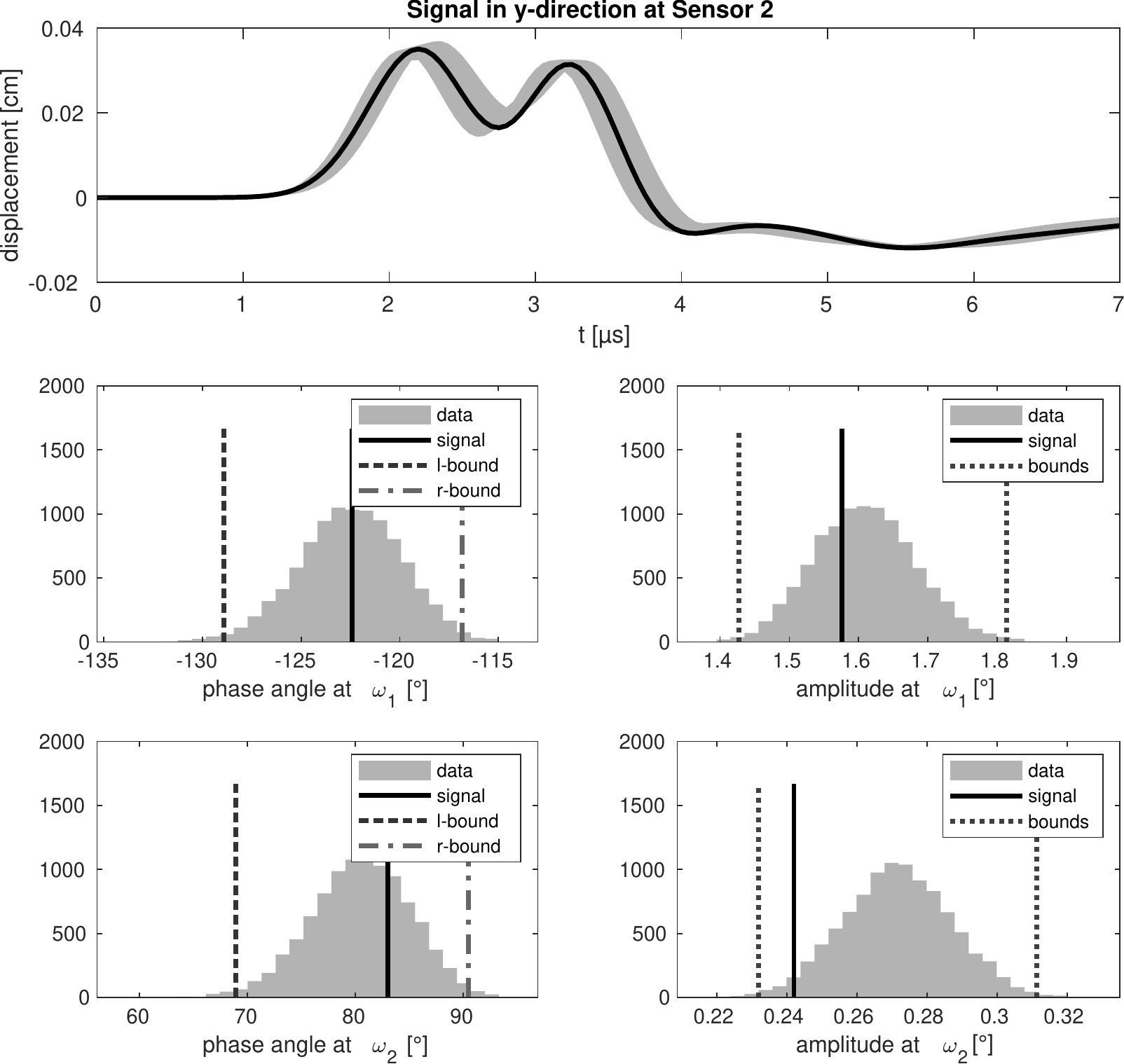}
	\caption{Scenario 1. Top row: Displacement in $y$-direction with Monte Carlo $99\%$-confidence band of solutions generated by Monte Carlo simulation. The measured solution (solid line) lies entirely within the band. Middle row: Histogram of corresponding features at frequency $\omega_1$, also showing bounds of the one-sided 99\% confidence regions (left) and the symmetric 99\% regions (right) as well as the value of the tested feature. The observed values lie well within the 99\% regions. Bottom row: Same as middle row, but for frequency $\omega_2$.}
	\label{fig:signal:fine}
\end{figure}

\begin{figure}[htb]\centering
	\noindent \hspace*{-4em}
	\begin{subfigure}{\linewidth+8em}
		\includegraphics[width=0.95\linewidth]{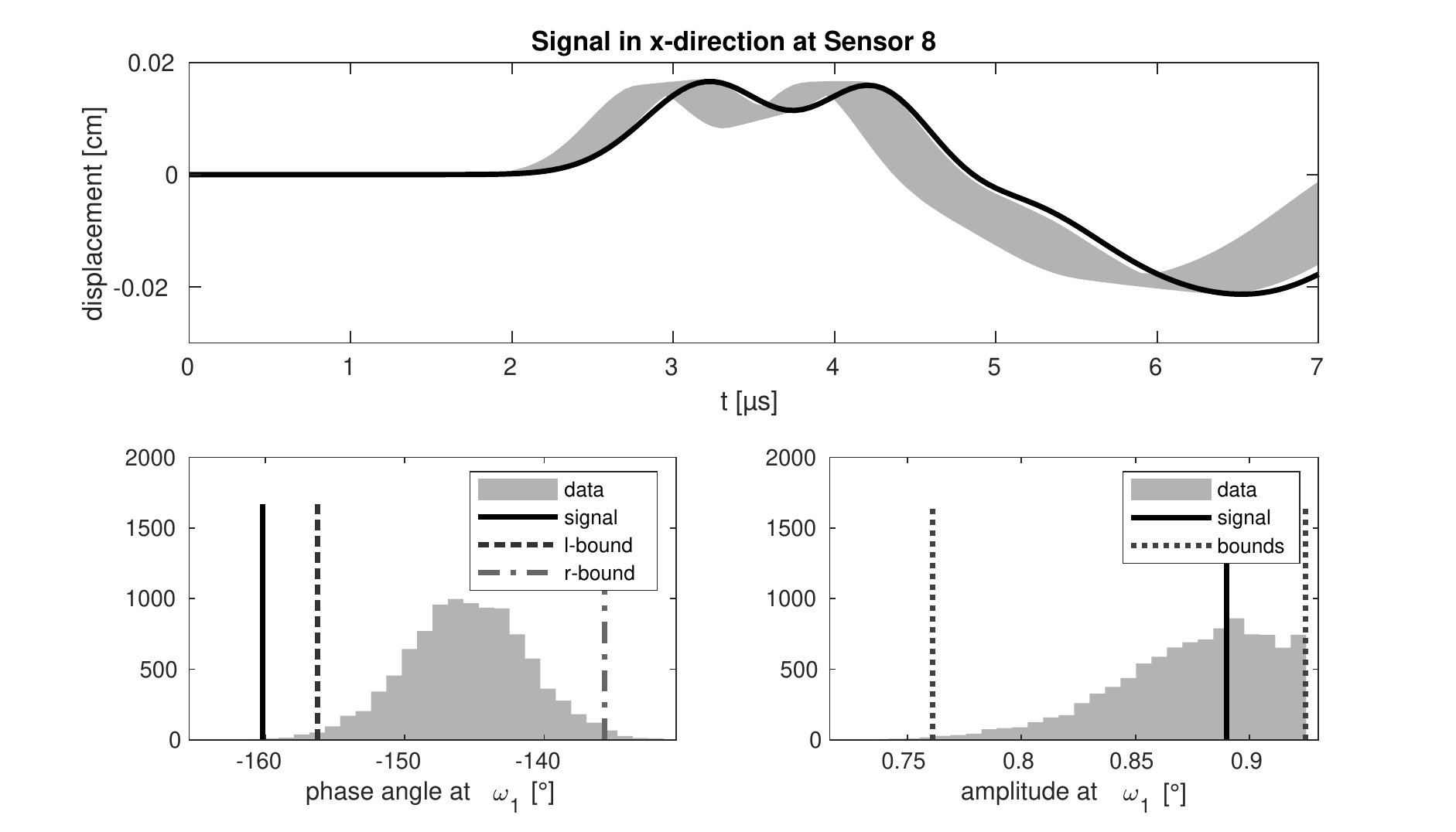}
\subcaption{}
	\end{subfigure}
	\noindent \hspace*{-4em}
	\begin{subfigure}{\linewidth+8em}
		\includegraphics[width=0.95\linewidth]{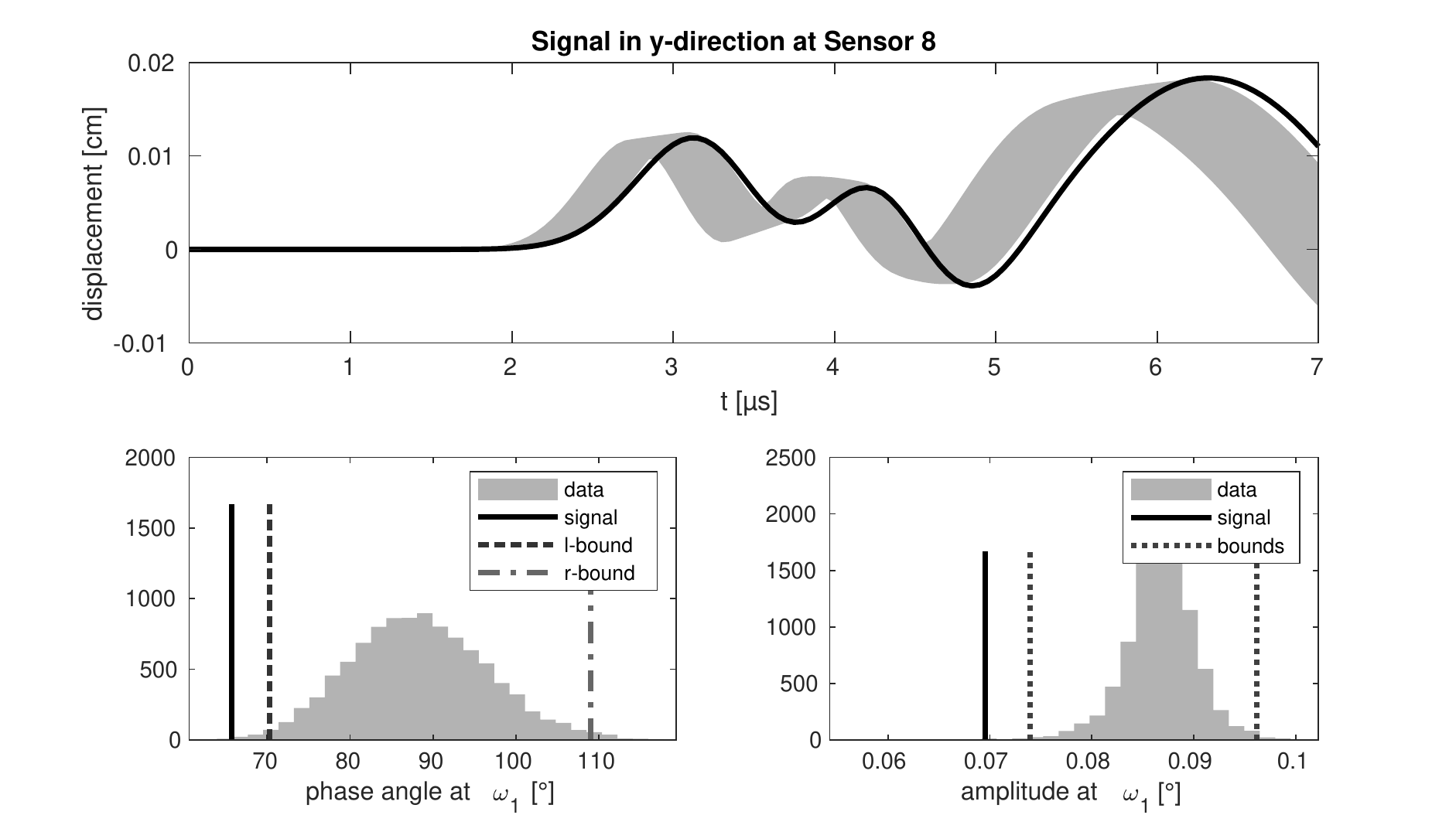}
\subcaption{}
	\end{subfigure}
	\caption{Scenario 2. (a) Displacement in $x$-direction with Monte Carlo $99\%$-confidence band of solutions and measured solution (solid line). Histogram of corresponding features at frequency $\omega_1$, also showing bounds of the one-sided 99\% confidence regions (left) and the symmetric 99\% regions (right) as well as the value of the tested feature.
(b) Same as part (a), but for displacement in $y$-direction. In both cases, the observed phase angle is smaller than the lower bound of the upper 99\% confidence region.}
	\label{fig:signal:smallE}
\end{figure}

\begin{figure}[htb]\centering
	\noindent \hspace*{-4em}
	\begin{subfigure}{\linewidth+8em}
		\includegraphics[width=0.95\linewidth]{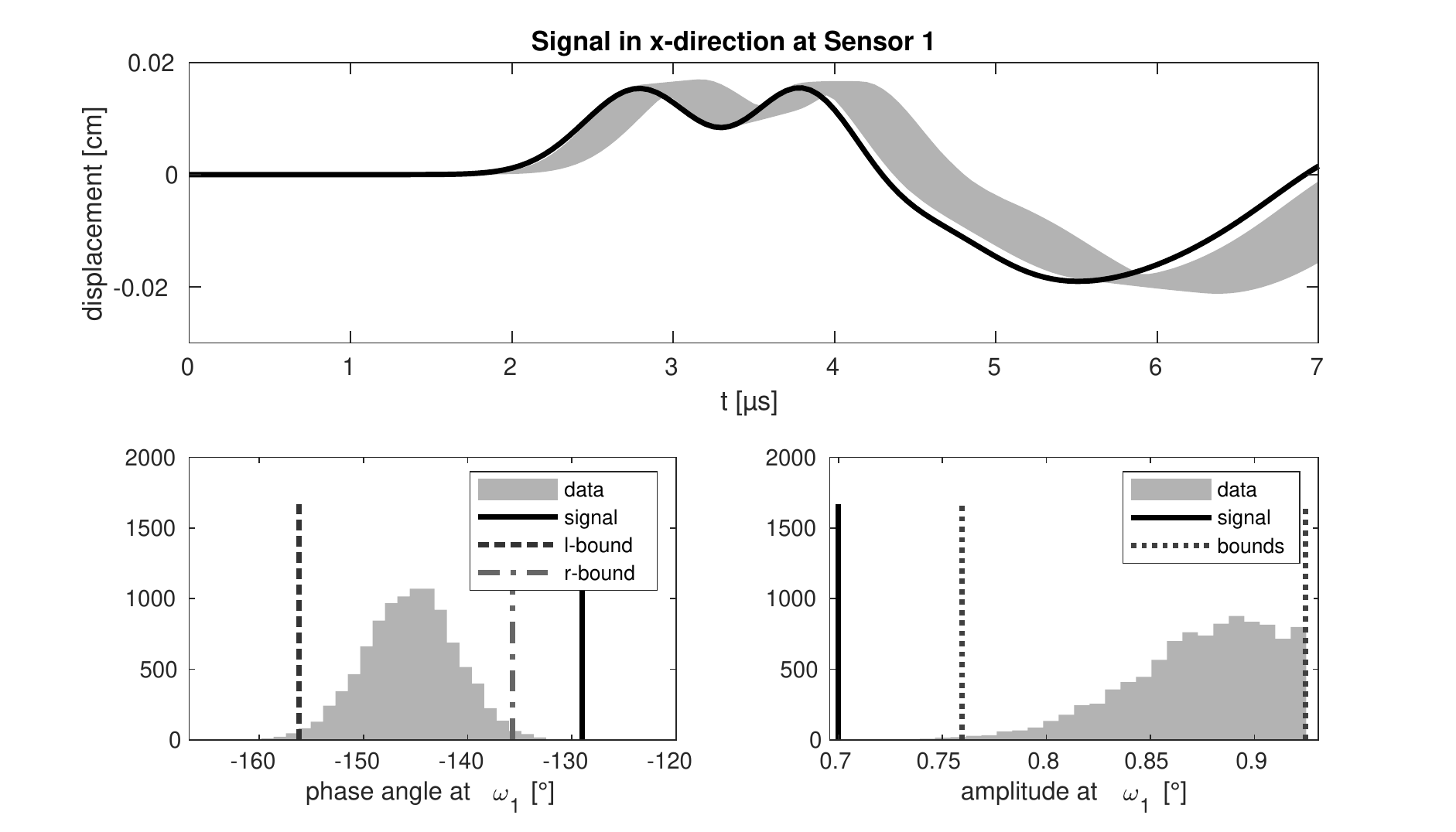}
\subcaption{}
	\end{subfigure}
	\noindent \hspace*{-4em}
	\begin{subfigure}{\linewidth+8em}
		\includegraphics[width=0.95\linewidth]{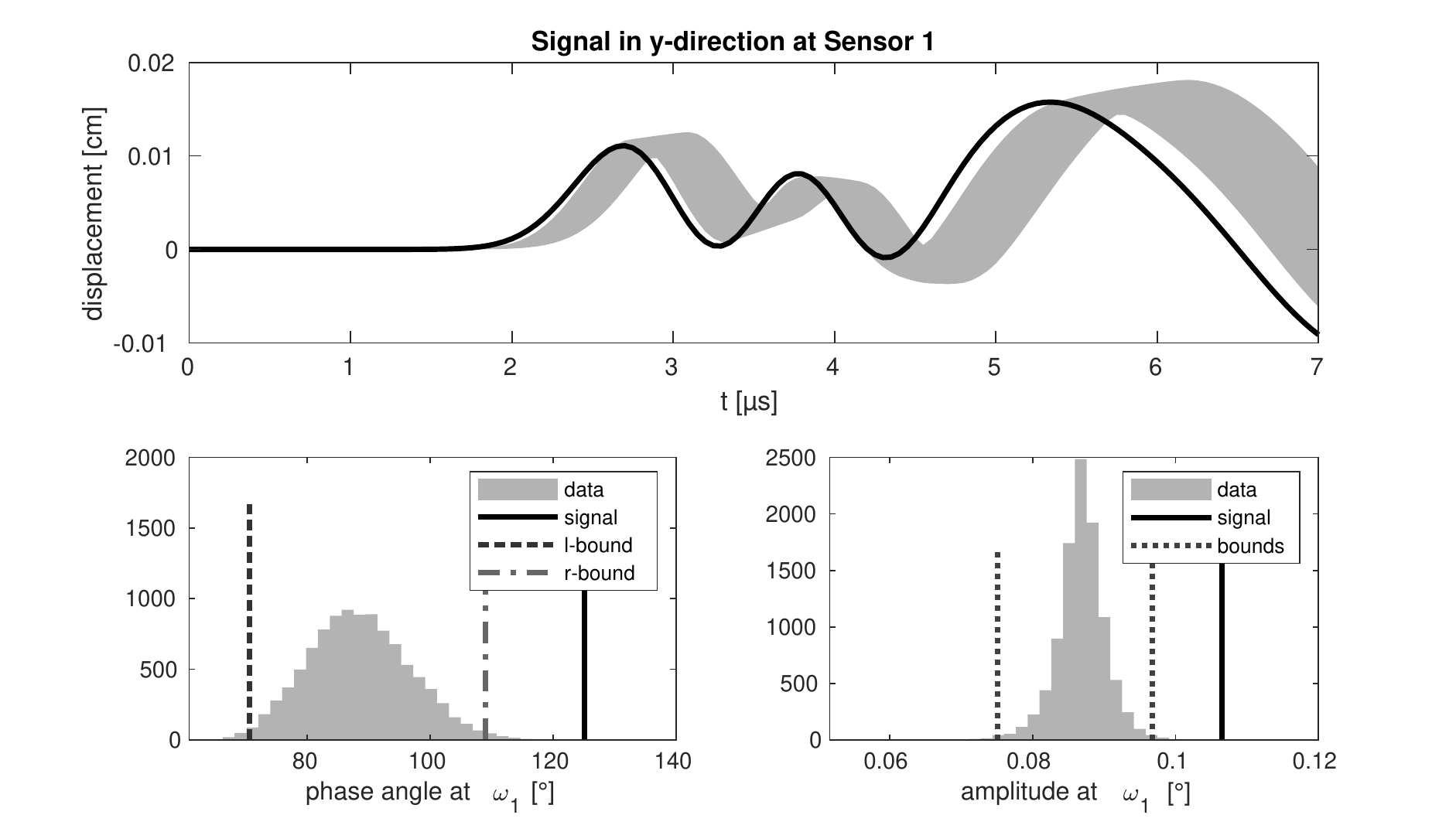}
\subcaption{}
	\end{subfigure}
	\caption{Scenario 3. (a) Displacement in $x$-direction with Monte Carlo $99\%$-confidence band of solutions and measured solution (solid line). Histogram of corresponding features at frequency $\omega_1$, also showing bounds of the one-sided 99\% confidence regions (left) and the symmetric 99\% regions (right) as well as the value of the tested feature.
(b) Same as part (a), but for displacement in $y$-direction. In both cases, the observed phase angle is larger than the upper bound of the lower 99\% confidence region.}
	\label{fig:signal:largeE}
\end{figure}

\begin{figure}[htb]\centering
	\noindent \hspace*{-4em}
	\begin{subfigure}{\linewidth+8em}
		\includegraphics[width=0.95\linewidth]{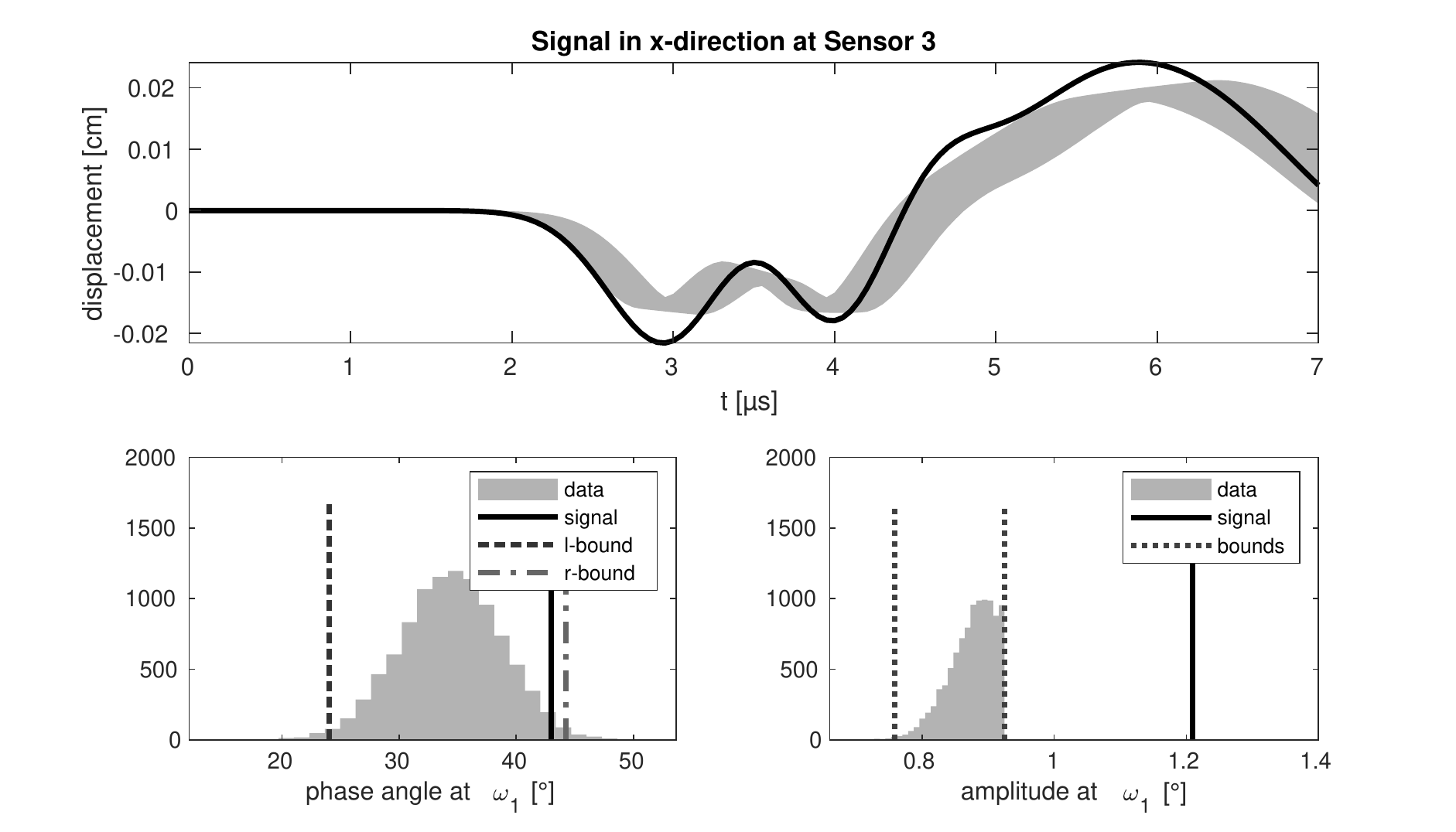}
\subcaption{}
	\end{subfigure}
	\noindent \hspace*{-4em}
	\begin{subfigure}{\linewidth+8em}
		\includegraphics[width=0.95\linewidth]{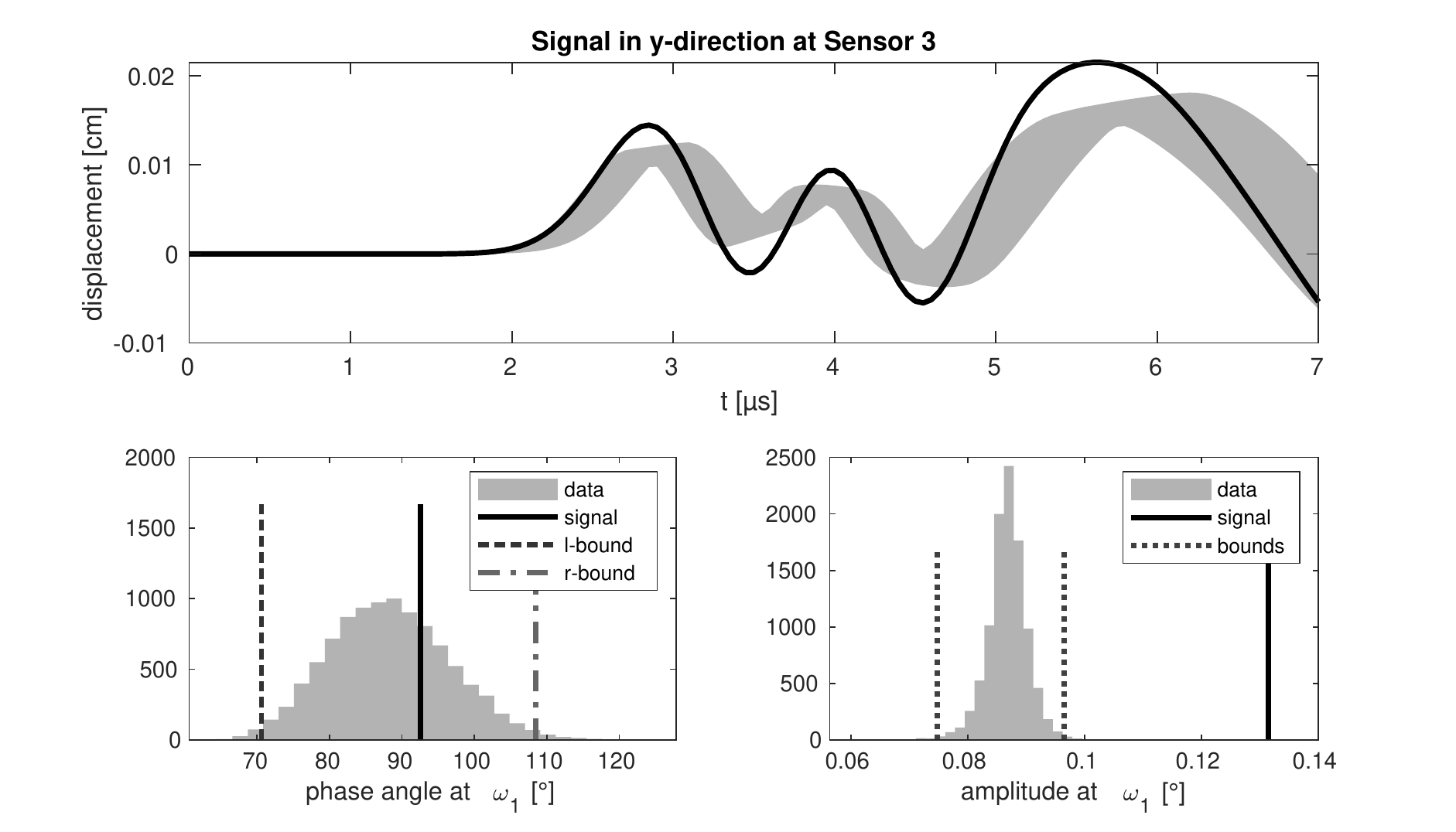}
\subcaption{}
	\end{subfigure}
	\caption{Scenario 4. (a) Displacement in $x$-direction with Monte Carlo $99\%$-confidence band of solutions and measured solution (solid line). Histogram of corresponding features at frequency $\omega_1$, also showing bounds of the one-sided 99\% confidence regions (left) and the symmetric 99\% regions (right) as well as the value of the tested feature.
(b) Same as part (a), but for displacement in $y$-direction. The crack strongly changes the signal in Sensor 3. The phase angles are inside the 99\% regions, but the amplitudes are by far out.}
	\label{fig:signal:crack}
\end{figure}

\begin{figure}[htb]
	\centering
	\begin{subfigure}{\linewidth}
		\includegraphics[width=0.88\linewidth]{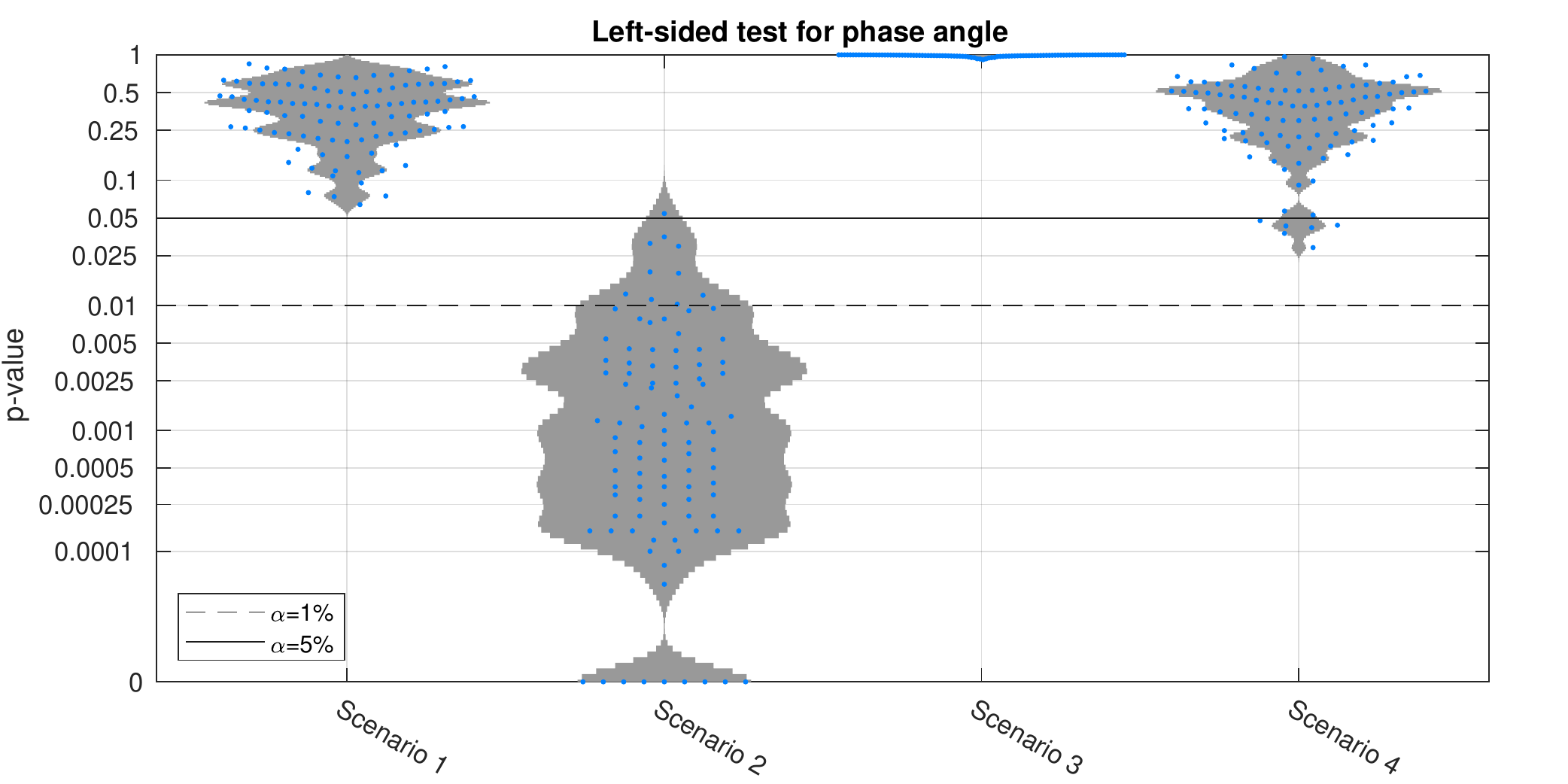}
		\subcaption{}
	\end{subfigure}
	\begin{subfigure}{\linewidth}
		\includegraphics[width=0.88\linewidth]{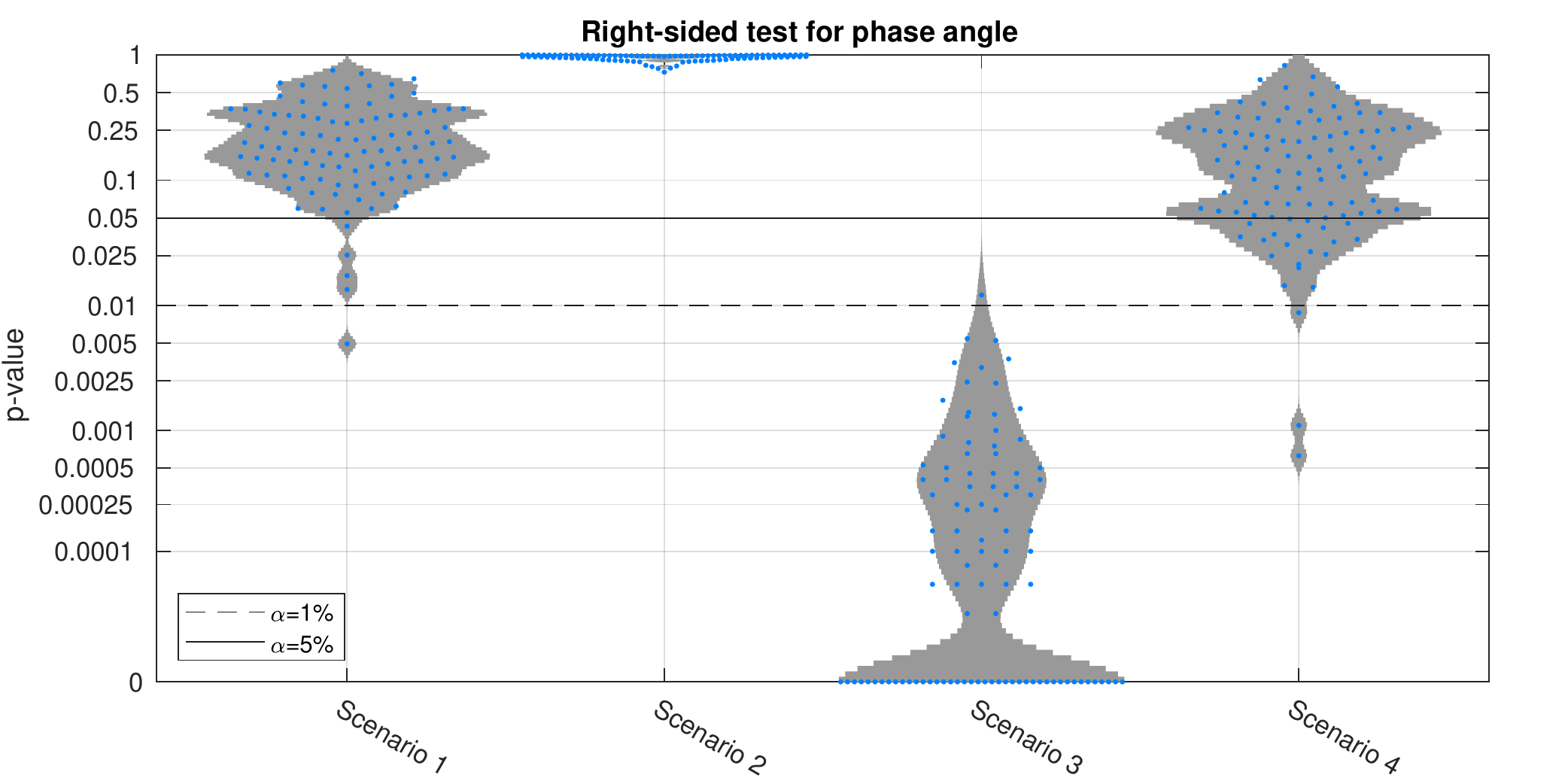}
		\subcaption{}
	\end{subfigure}
	\begin{subfigure}{\linewidth}
		\includegraphics[width=0.88\linewidth]{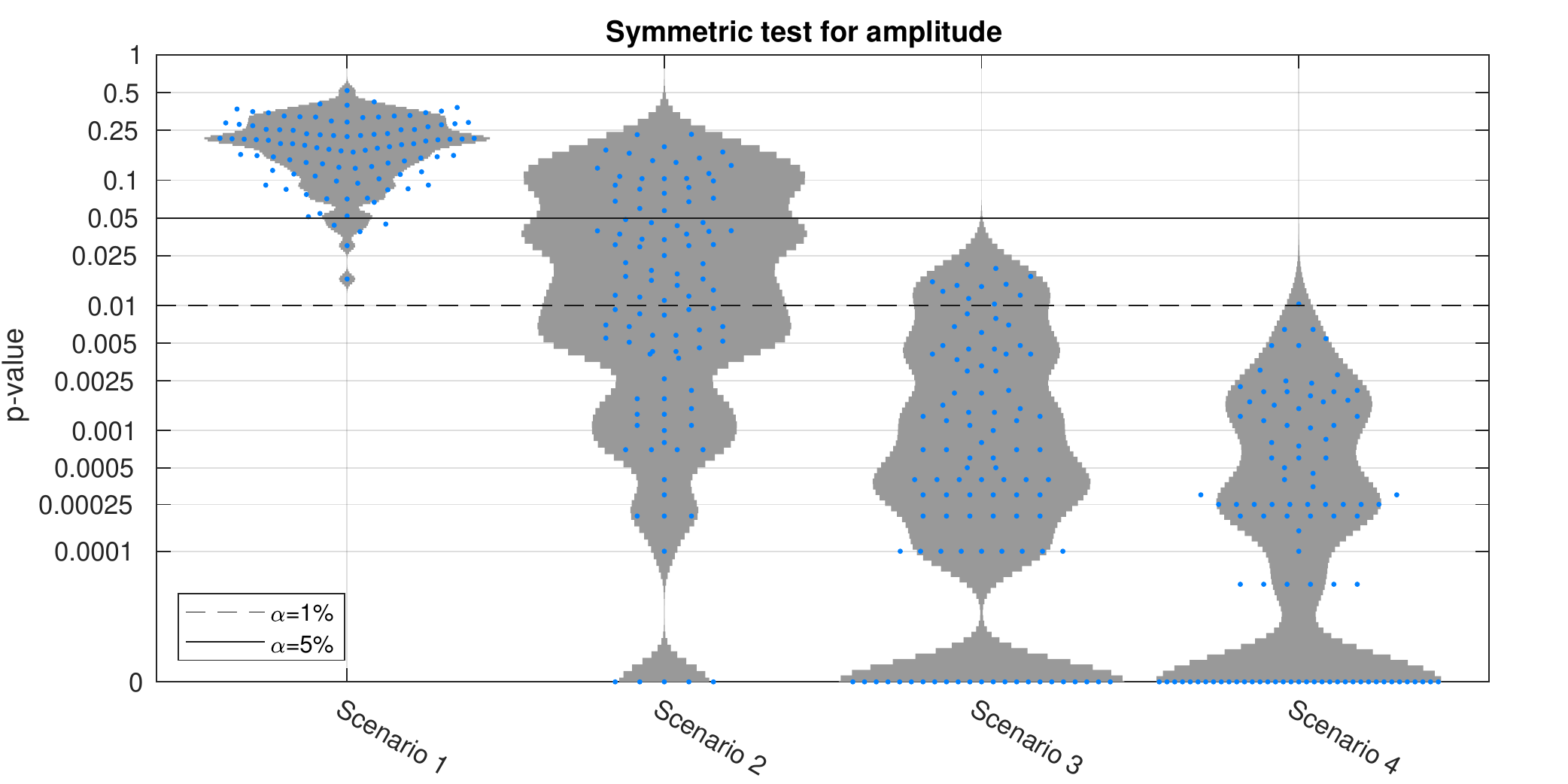}
		\subcaption{}
	\end{subfigure}
	\caption{Violin plots of the computed $p$-values. The dots represent the $p$-values of a single FE-simulation. The gray shaded regions indicate an estimate of the distribution of the $p$-values. The horizontal thin and dashed lines represent the 5\% and 1\%-level, respectively. (a) Test (I); (b) Test (II); (c) Test (III).}
	\label{fig:tests_scens}
\end{figure}

\begin{figure}[htb]
	\centering
	\begin{subfigure}{\linewidth}
		\includegraphics[width=0.88\linewidth]{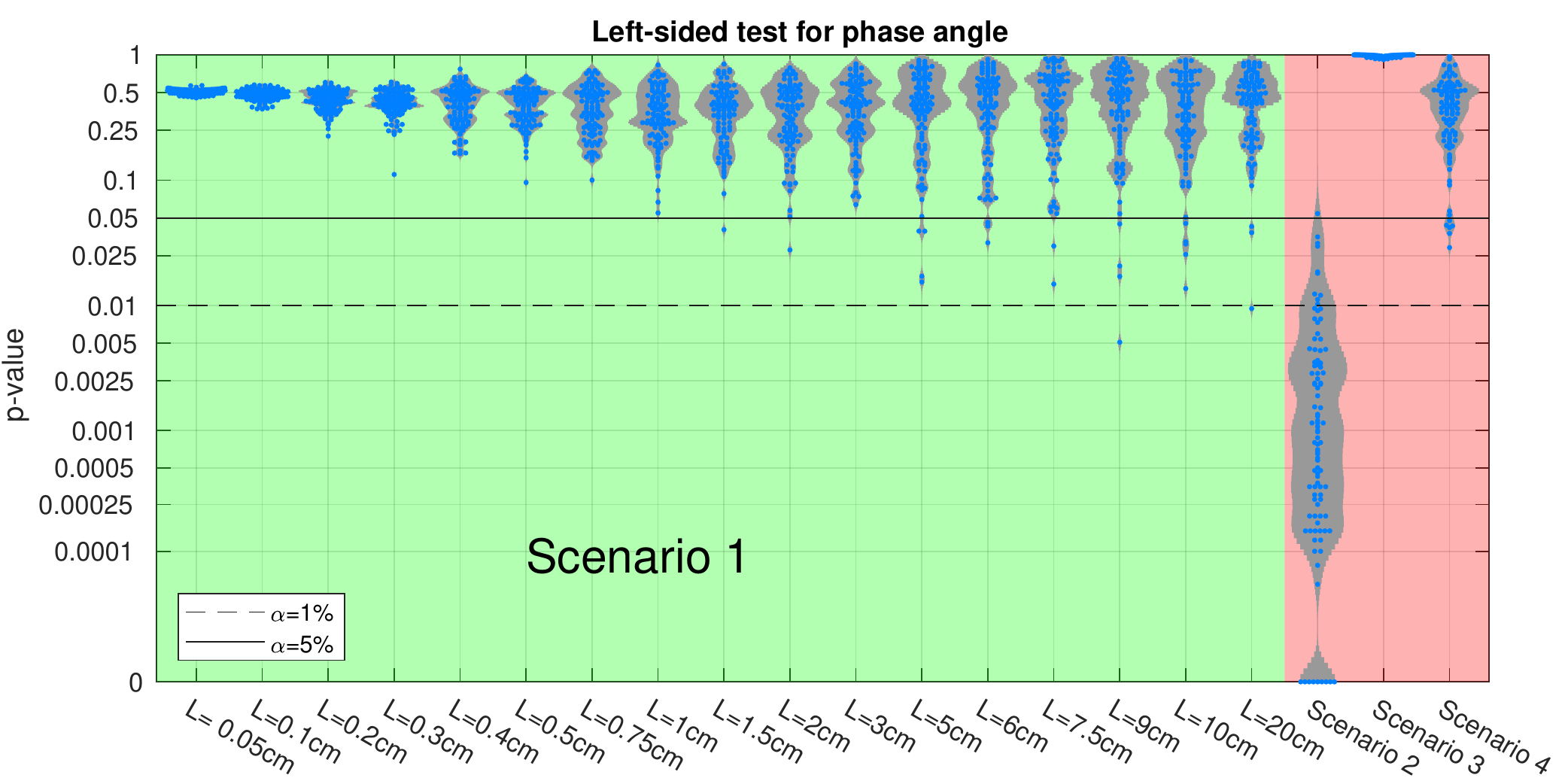}
		\subcaption{}
	\end{subfigure}
	\begin{subfigure}{\linewidth}
		\includegraphics[width=0.88\linewidth]{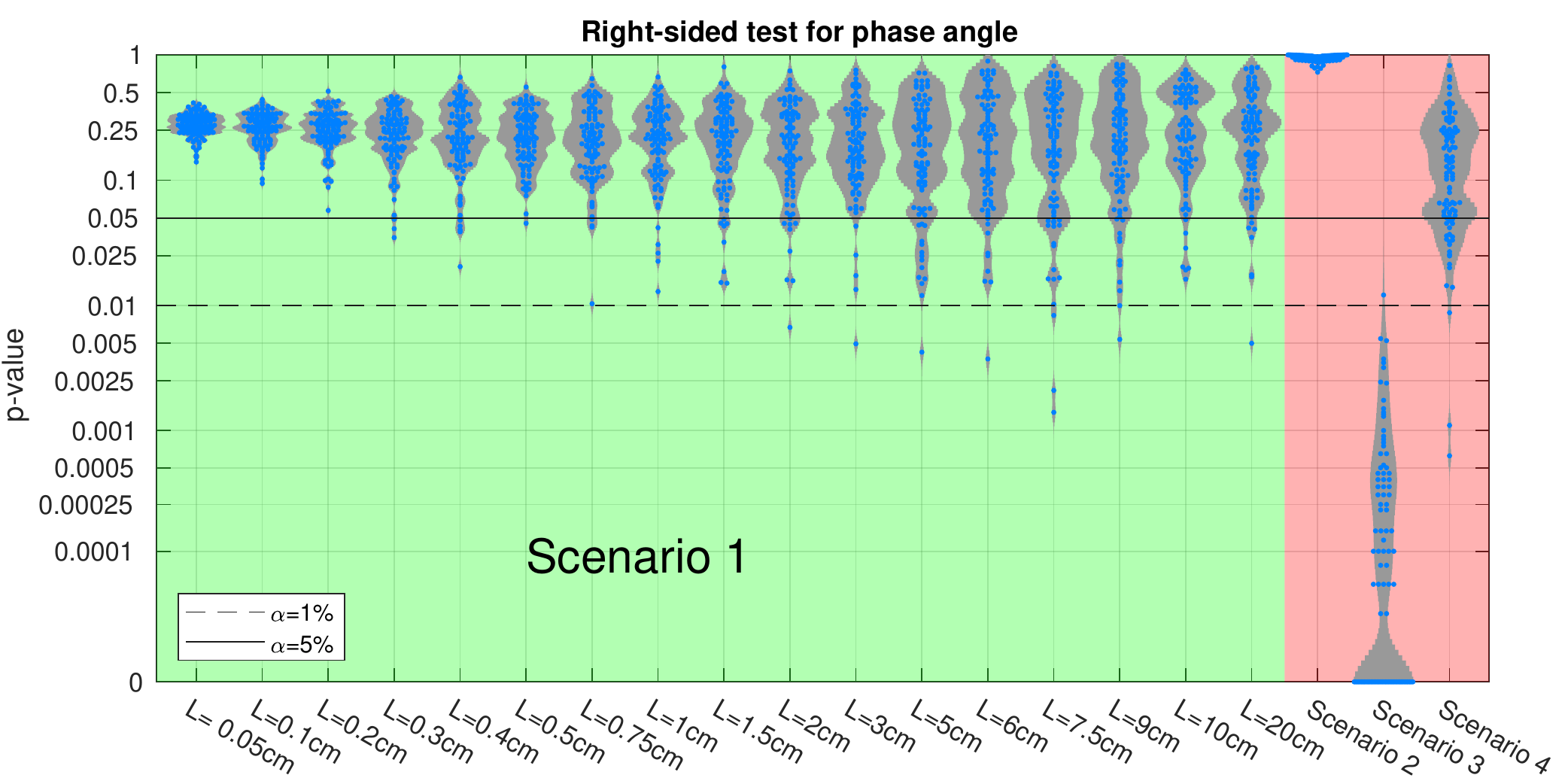}
		\subcaption{}
	\end{subfigure}
	\begin{subfigure}{\linewidth}
		\includegraphics[width=0.88\linewidth]{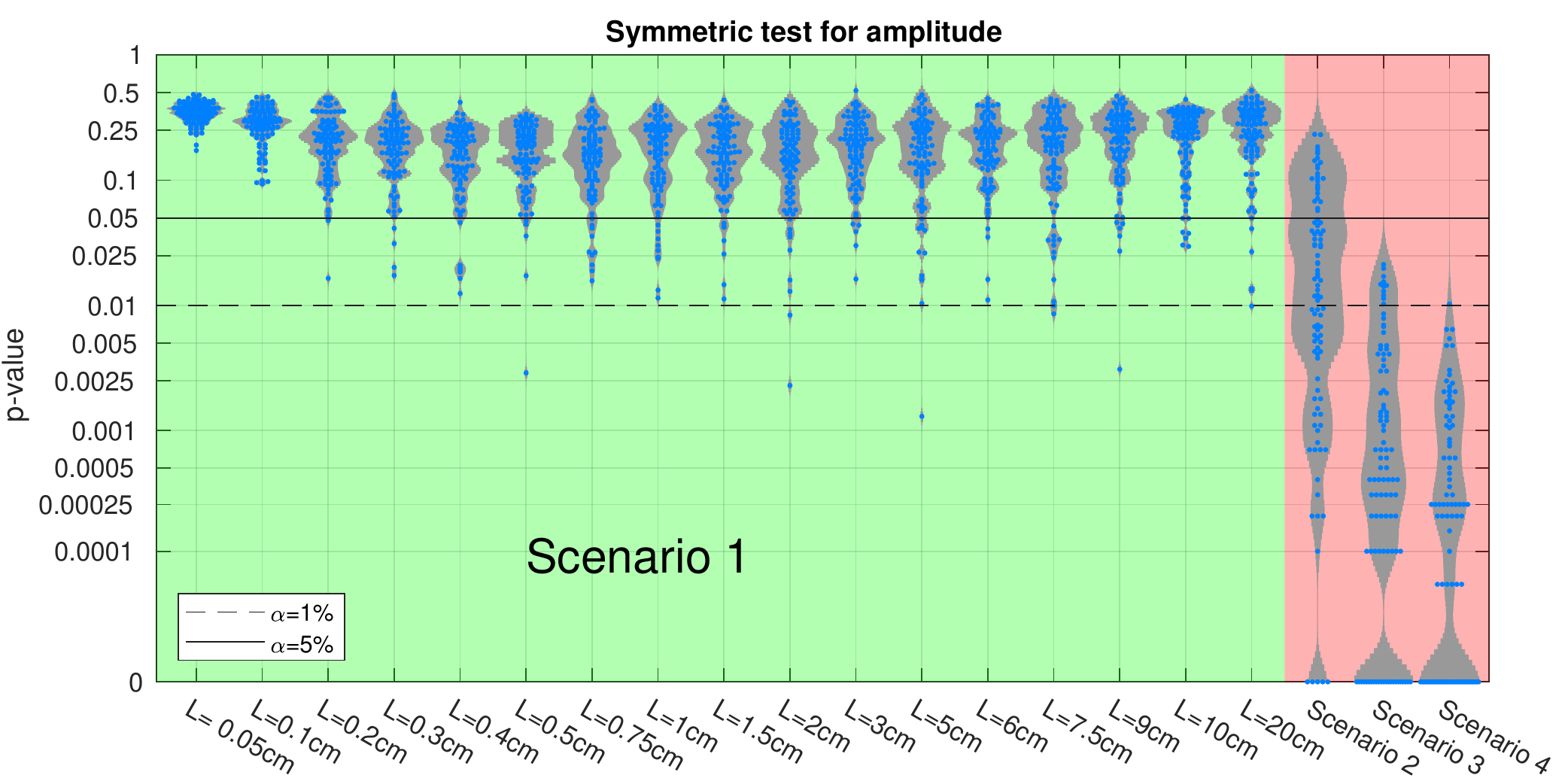}
		\subcaption{}
	\end{subfigure}
	\caption{Violin plots of the computed $p$-values, showing dependence on correlation length of undamaged material. The dots represent the $p$-values of a single FE-simulation. The gray shaded regions indicate an estimate of the distribution of the $p$-values. The horizontal thin and dashed lines represent the 5\% and 1\%-level, respectively. The highlighted areas to the right of each rectangle recall the results for Scenarios 2 -- 4 with $L_E = 3$ cm. (a) Test (I); (b) Test (II); (c) Test (III).}
	\label{fig:tests_Ls}
\end{figure}

\begin{figure}[htb]
	\centering
	\begin{subfigure}{\linewidth}
		\includegraphics[width=0.88\linewidth]{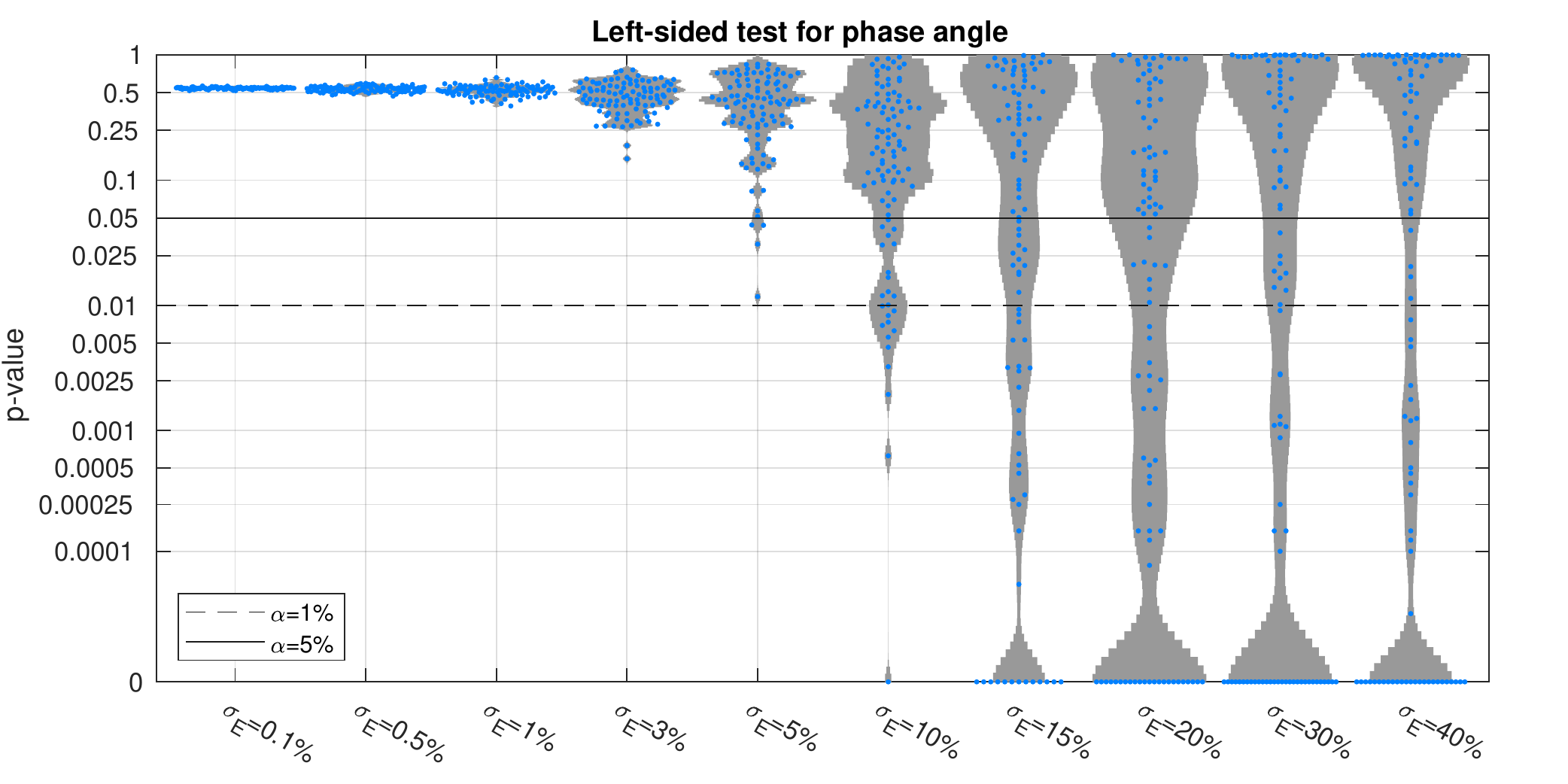}
		\subcaption{}
	\end{subfigure}
	\begin{subfigure}{\linewidth}
		\includegraphics[width=0.88\linewidth]{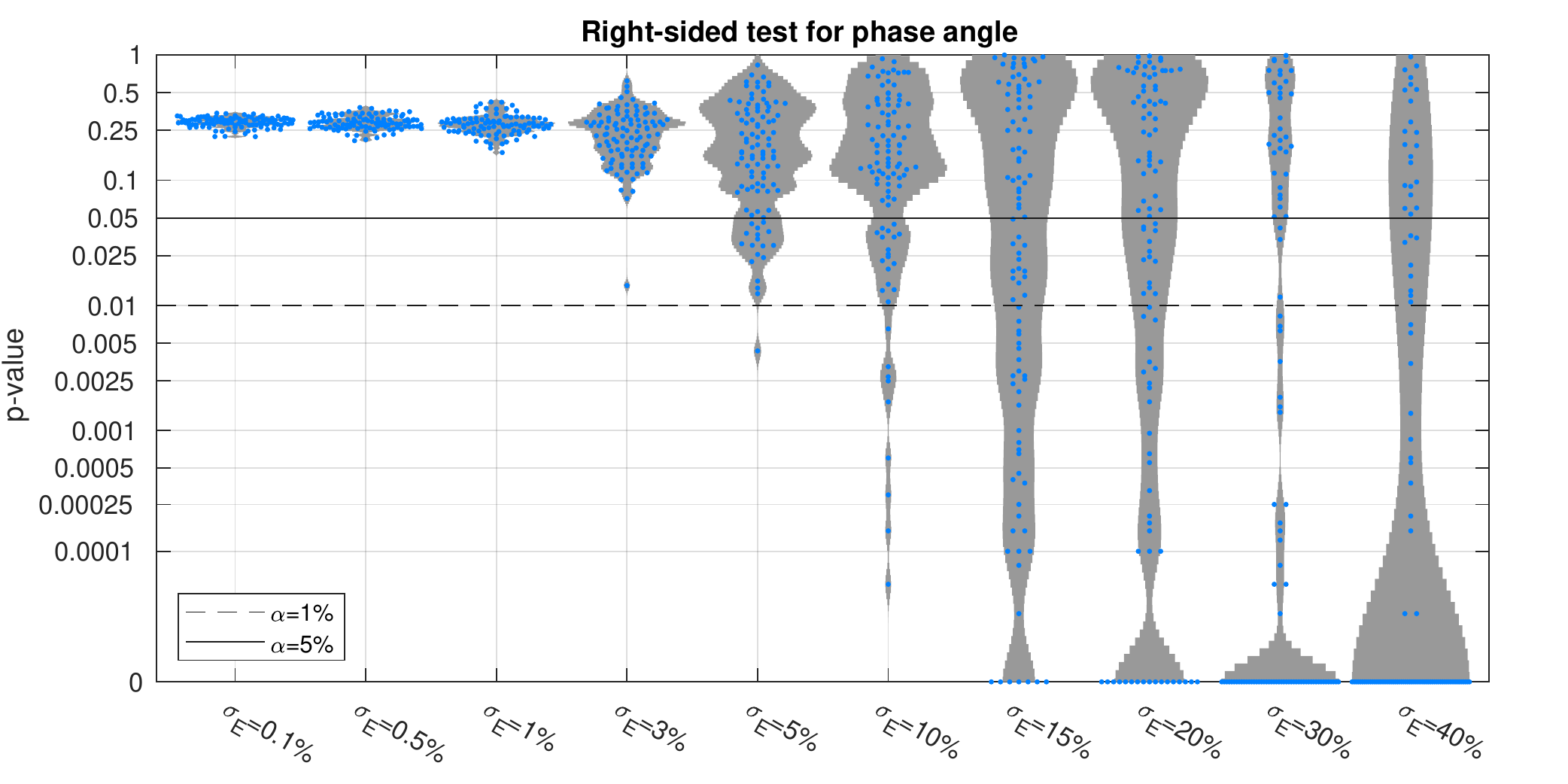}
		\subcaption{}
	\end{subfigure}
	\begin{subfigure}{\linewidth}
		\includegraphics[width=0.88\linewidth]{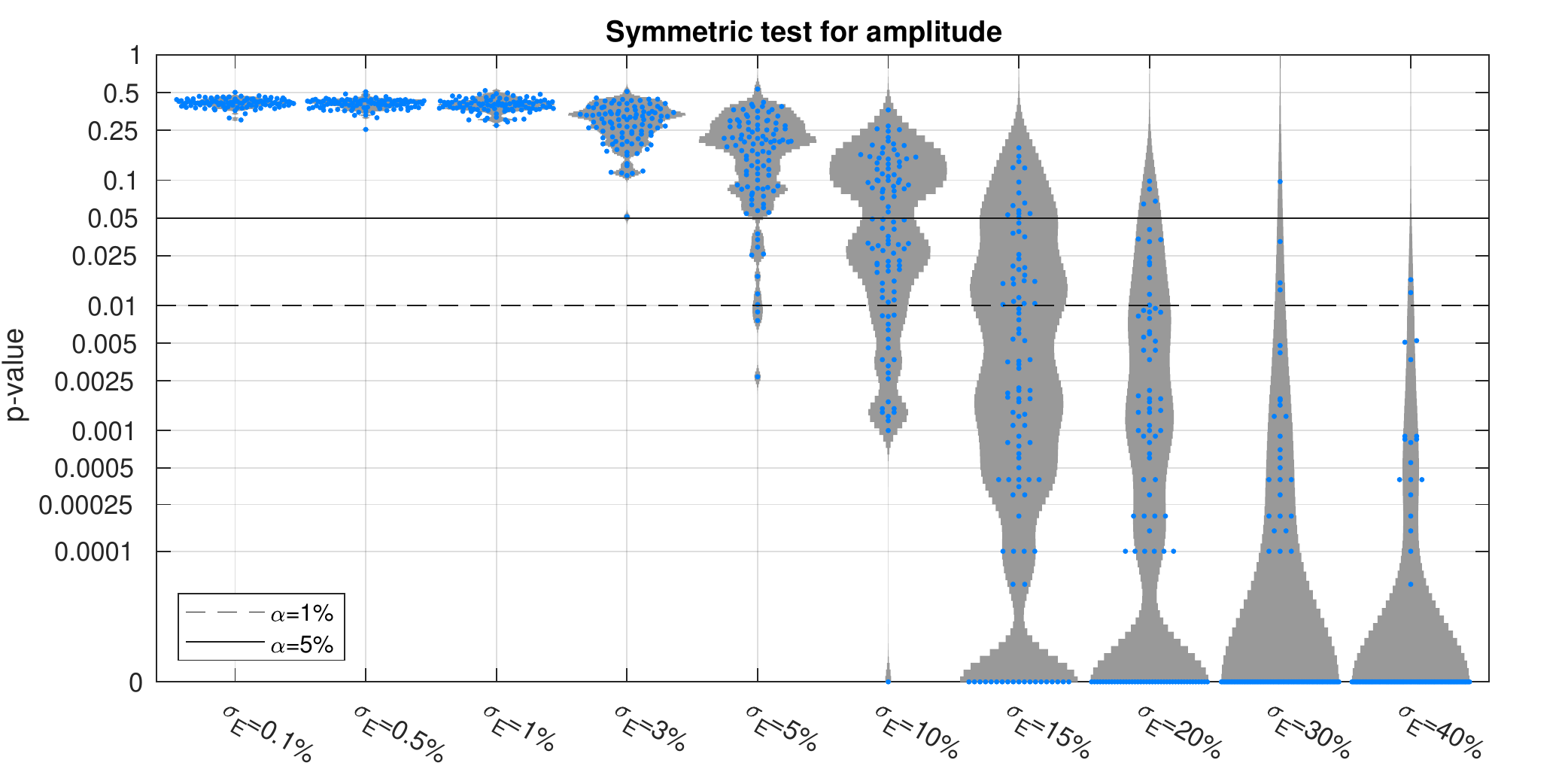}
		\subcaption{}
	\end{subfigure}
	\caption{Violin plots of the computed $p$-values, showing dependence on standard deviation of undamaged material. The dots represent the $p$-values of a single FE-simulation. The gray shaded regions indicate an estimate of the distribution of the $p$-values. (a) Test (I); (b) Test (II); (c) Test (III).}
	\label{fig:test_sigma}
\end{figure}

\begin{landscape}
	\begin{table}[htb]
		\centering
		\begin{tabular} {|c|c||c|c|c|c|c|c|c|c|c|c|c|c|c|c|c|c|c|c|}\hline
			\multicolumn{2}{|c||}{} & \multicolumn{17}{c|}{$L_E$ [cm]}\\\hline
			Test&$\alpha$& 0.05&0.1&0.2&0.3&0.4&0.5&0.75&1&1.5&2&3&5&6&7.5&9&10&20\\\hline	\hline
			(I)&5\%& 0  & 0  & 0  & 0  & 0  & 0  & 0  & 0  & 1  & 1  & 0  & 4  & 4  & 2  & 4  & 5  & 3\\\hline
			(I)& 1\%& 0  & 0  & 0  & 0  & 0  & 0  & 0  & 0  & 0  & 0  & 0  & 0  & 0  & 0  & 1  & 0  & 1\\\hline
			(II)&5\%& 0  & 0  & 0  & 0  & 0  & 0  & 0  & 0  & 1  & 1  & 0  & 4  & 4  & 2  & 4  & 5  & 3 \\\hline
			(II)& 1\%& 0  & 0  & 0  & 0  & 0  & 0  & 0  & 0  & 0  & 0  & 0  & 0  & 0  & 0  & 1  & 0  & 1\\\hline
			(III)&5\%& 0  & 0  & 0  & 3  & 5  & 1  & 4  & 5  & 8  & 9  & 5  & 14  & 9  & 17  & 11  & 7  & 7 \\\hline
			(III)& 1\%& 0  & 0  & 0  & 0  & 0  & 0  & 0  & 0  & 0  & 1  & 1  & 1  & 1  & 3  & 1  & 0  & 1   \\\hline
		\end{tabular}
		\caption{Number of $p$-values below the $\alpha$-level among 100 repetitions of the tests (false rejections of null hypothesis), depending on the correlation length of $E$.}\vspace{4em}
		\label{tbl:pvalues_of_L}
		\centering
		\begin{tabular} {|c|c||c|c|c|c|c|c|c|c|c|c|c|c|c|c|c|c|c|c|}\hline
			\multicolumn{2}{|c||}{} & \multicolumn{10}{c|}{Coefficient of variation of $E$}\\\hline
			Test&$\alpha$&0.1\%&   0.5\%&    1\%&    3\%&    5\%&  10\%&  15\%&   20\%&   30\%&   40\% \\\hline	\hline
			(I)&5\%   & 0  & 0  & 0  & 0  & 4  & 24  & 45  & 52  & 51  & 46   \\\hline
			(I)& 1\%  & 0  & 0  & 0  & 0  & 0  & 12  & 33  & 43  & 42  & 42   \\\hline
			(II)&5\%  & 0  & 0  & 0  & 1  & 18  & 26  & 51  & 50  & 67  & 78   \\\hline
			(II)& 1\%  & 0  & 0  & 0  & 0  & 1  & 9  & 37  & 37  & 64  & 70   \\\hline
			(III)&5\% & 0  & 0  & 0  & 0  & 11  & 54  & 86  & 96  & 99  & 100   \\\hline
			(III)& 1\% & 0  & 0  & 0  & 0  & 3  & 21  & 66  & 86  & 96  & 98   \\\hline
		\end{tabular}
		\caption{Number of $p$-values below the $\alpha$-level among 100 repetitions of the tests (false rejections of null hypothesis), depending on the standard deviation of $E$.}
		\label{tbl:pvalues_of_sigma}
	\end{table}
\end{landscape}

\end{appendix}



\end{document}